\documentclass[a4paper,12pt]{article}
\usepackage[a4paper,vmargin={3.5cm,3.5cm},hmargin={2.5cm,2.5cm}]{geometry}
\usepackage[font=sf, labelfont={sf,bf}, margin=1cm]{caption}
\usepackage{amsmath,amssymb,amsfonts,amsthm}
\usepackage{graphicx}

\newtheorem{corollary}{Corollary}
\newtheorem{theorem}[corollary]{Theorem}

\newtheorem{lemma}[corollary]{Lemma}

\newtheorem{proposition}[corollary]{Proposition}

\newtheorem{remark}{Remark}

\newcommand{\M}{\mathcal{M}}
\newcommand{\R}{\mathbb{R}}

\newcommand{\CC}{\mathcal{C}}
\newcommand{\D}{\mathbb{D}}

\newcommand{\N}{\mathbb{N}}
\renewcommand{\H}{\mathbb{H}}

\newcommand{\A}{\mathcal{A}}
\newcommand{\T}{\mathcal{T}}
\newcommand{\TT}{\mathbf{T}}
\newcommand{\J}{\mathcal{J}}

\newcommand{\op}[1]{\operatorname{#1 }}

\renewcommand{\P}{\mathcal{P}}
\renewcommand{\i}{ \mathrm{i}}

\newcommand{\HH}{\mathbb{H}}
\newcommand{\xii}{\chi}

\title{The Markovian hyperbolic triangulation}
\author {Nicolas Curien \and Wendelin Werner}
\date {Ecole Normale Sup\'erieure and Universit\'e Paris-Sud}

\begin{document}
\maketitle
\begin {abstract}
 We construct and study the unique random tiling of the hyperbolic plane into ideal hyperbolic triangles (with the three corners located on the boundary) that is invariant (in law) with respect to M\"obius transformations, and possesses a natural spatial Markov property that can be  roughly described as the conditional independence of the two parts of the triangulation on the two sides of the edge of one of its triangles. 
\end {abstract}

\section {Introduction}

The study of the scaling limit of critical two-dimensional discrete models from statistical physics has given rise to various 
random objects in the continuum that combine conformal invariance with a ``spatial Markov property'' that is inherited from the locality of the interactions 
 in the discrete models (one can think of course about Schramm's SLE processes \cite{S00}). 

In the present paper we shall exhibit and study a special M\"obius-invariant random triangulation of the Poincar\'e disk $\D$ endowed with its hyperbolic complex structure, that
 possesses a certain  spatial Markov property. 
Let us first very briefly explain what type of triangulations we have in mind. A (hyperbolic) triangle $T$ will be determined by its three corners, that we will always take on $\partial \D$, and $T$ will be the ``inside'' of the three hyperbolic lines joining these three points (recall that these hyperbolic lines are circular chords when viewed in the Euclidean setting).   
We say that $\mathbf {T}$ is a complete hyperbolic triangulation of $\D$ if it is a disjoint collection of such triangles, and if the complement of the union of all these triangles has zero hyperbolic measure. 

We say that a random triangulation $\mathbf {T}$ is M\"obius-invariant if its law is invariant under all conformal transformations from $\D$ onto itself. In other words, for any M\"obius transformation $\phi$ of the unit disk onto itself, the law of $\phi ( \mathbf {T})$ is the same as that of $\mathbf {T}$.
Our main statement can be described as follows:

\medbreak

{\sl There exists a unique random complete M\"obius-invariant triangulation of $\mathbb{D}$ that fulfills a spatial Markov property that can be loosely speaking described as follows: Given a triangle $T=(abc)$ in this triangulation $\TT$ (in fact the rigorous statement is to say that $T$ is the triangle that contains the origin in $\TT$), the triangulation $\TT$ restricted to the three connected components of the complement of $T$ in $ \mathbb{D}$ are conditionally independent, and moreover, the part that is beyond $(bc)$ is independent of the position of $a$.}

\begin{figure}[!ht]
 \begin{center}
  \includegraphics[height=8cm]{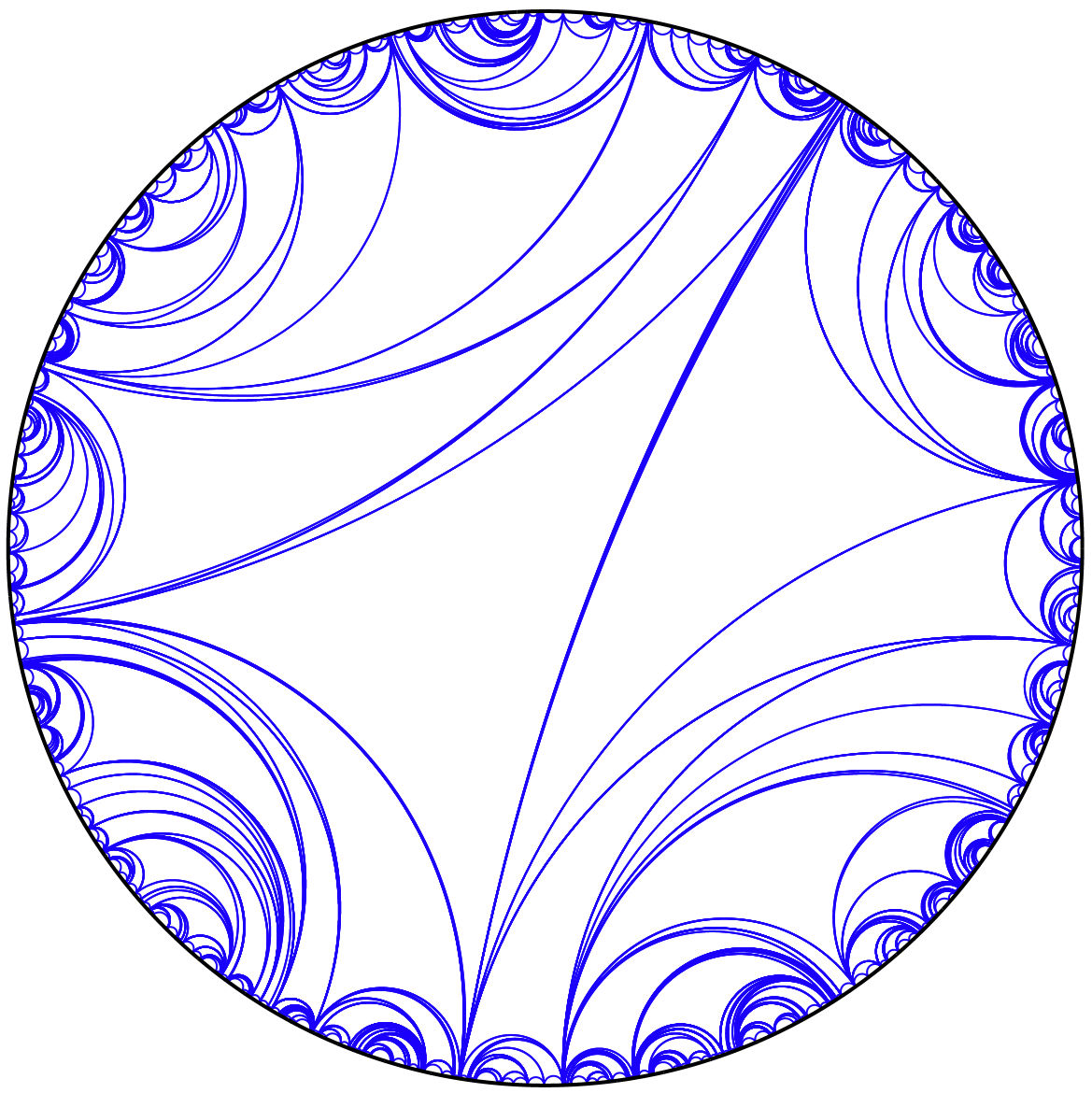}
 \caption {\label {diskpic}Sample of our triangulation in the disk.}
 \end{center}
 \end{figure}

\medbreak
This will be stated more rigorously in Theorem \ref {main}. Uniqueness means of course here uniqueness of the law of the triangulation. Heuristically, the spatial Markov property means that there is conditional independence of both sides of an edge in the triangulation, so that the role of the edges in our triangulations is reminiscent of that of an interface in a nearest-neighbor interaction model from statistical physics.

\medbreak

Discrete models, such as triangulations of convex polygons have been thoroughly studied in combinatorics, physics or geometry. Some triangulations of the disk can be viewed as continuous counterparts to these discrete models, and various random triangulations of the disk have been defined and studied, in particular in recent years (see for instance \cite{Ald94b,CLG10} and the references therein). 
However, the particular random triangulation that we construct and study in the present paper is different (we would like to stress that it is {\em not} the same as the uniform triangulation defined by Aldous \cite{Ald94b} that can be viewed as the scaling limit of the uniform triangulation of a $n$-gon; we shall for instance see that our triangulation is much thinner), and despite its rather striking properties it does (to our knowledge) not seem to have been studied before.

\medbreak
In order to help getting a feeling of what is going on, let us provide a brief heuristic discussion. Assume that a random triangulation $\TT$ is complete, M\"obius-invariant and Markovian.  We can first sample the triangle $T(0)$ of $\TT$ that contains the origin -- and it will be easy  (see Section \ref {S2}) to identify its law from the  conformal invariance and completeness of $\TT$. We write $a,b$ and $c$ for the three apexes of $T(0)$ ordered anti-clockwise. Then, we can start exploring the three pieces of the complement of $T(0)$ independently, because of the spatial Markov property. A first naive guess is that the edge $(bc)$ will also be  one of the edges of another triangle of $\TT$, that ``neighbors'' $T(0)$. We can wonder what the conditional law of its third corner $a'$ will be. Conformal invariance (and the fact that $a'$ is conditionally independent of $a$ given $(bc)$) imposes that this (conditional) law is invariant under all M\"obius transformations of the disc that fix $b$ and $c$. But all (non-zero) measures supported on the arc of $\partial \mathbb{D}$ between $b$ and $c$ that are invariant under all these transformations  necessarily have an infinite mass (in fact, they are the  multiples of the image of the measure  $\xi(\mathrm {d} x) := \mathrm{d}x/x$ on $\R_{+}$ under any fixed M\"obius transformation that maps the upper half plane  onto the unit disc, and $0$ and $\infty$ onto $b$ and $c$ respectively), and more precisely this infinite mass lies in the neighborhood of $b$ and of $c$. So, this attempt to construct a neighboring triangle  in a M\"obius-invariant way fails, but it suggests to those of us who are acquainted to L\'evy processes which way to go: When exploring 
the triangles ``outwards'' starting from the $(bc)$ boundary of $T(0)$, one will use a Poisson point process, with intensity given by $\xi$, that will be used at each ``time'' to choose the following new corners. In particular (because $\xi$ is an infinite measure), almost surely, two different triangles $T$ and $T'$ in the triangulation will never be adjacent (there will always be infinitely many other very thin -- in the Euclidean meaning -- triangles that are separating $T$ from $T'$). In fact, it will turn out that there are not that many triangles either: For any large $n$, the number of triangles of (Euclidean) width in $[2^{-n-1}, 2^{-n})$ that are separating $T$ and $T'$ is of constant order (it is random but its mean is roughly constant). This explains why one could at first glance think that big triangles can happen to be adjacent to each other by looking at the simulation depicted in Figure \ref {diskpic}. 

\medbreak
The paper is organized as follows. In Section \ref {S2}, we collect and derive some rather general or elementary facts, we write down definitions and state our main result, Theorem \ref {main}.
In Section \ref {S3}, inspired by the previous heuristic, we show that if $\TT$ is a M\"obius-invariant complete and Markovian triangulation, then it necessarily corresponds to some Poisson point process that we describe. This argument will prove the fact that the law of $\TT$ is {\em unique}, if it exists. In Section \ref {S4}, we define explicitly a random triangulation (again, using a Poisson point process), and we check that it indeed satisfies all the required properties, so that {\em existence} of {the} M\"obius-invariant complete Markovian triangulation follows.  

In Section \ref {S5}, we discuss and state results dealing with M\"obius-invariant Markovian tilings of the disk into other hyperbolic polygons than triangles (analogous statements and constructions exist for instance for tilings into conformal squares). 
In the final section, we make a few comments, and list a couple of open questions. 

We are going to assume that the reader is acquainted with basic properties of M\"obius transformations on the one hand, and basic knowledge about Poisson point processes, pure jump processes and subordinators (as can be found in \cite {Ber96,Ber99,RY99}) on the other hand. As we are aware that this is not such an usual mix of backgrounds, we will however try to recall some of the basic features that we will use.

\section {Simple preliminaries}
\label {S2}

\subsection {Hyperbolic triangles}
\label{preliminaires}

We will mostly use the unit disc $\D = \{ z \in \mathbb{C} : |z|<1\}$ to represent the hyperbolic plane. At some point,
 it will also be convenient to work in the upper half-plane $\H = \{ z \in \mathbb{C} : \ \mathrm{Im}(z)>0\}$. Throughout the paper, $\psi^{-1}$ will denote the conformal map from $\H$ onto $\D$ defined by  
$\psi^{-1} (z) = (z-\i) / (z+\i)$, which maps $\i$ onto the origin and infinity onto $1$. 

For any pair of \emph{distinct} points $a$ and $b$ on $\partial \D$, we define the hyperbolic line $(ab)$ in $\D$ to be the circular chord in $\overline \D$ that crosses $\partial \D$ orthogonally at both $a$ and $b$ (when $a=-b$, this  ``circular chord'' is in fact a diameter line). In order to try to avoid confusions, we will use $[a,b]$ to denote straight Euclidean segments. 

If we consider three distinct points $a$, $b$ and $c$ on $\partial \D$, we can define a hyperbolic triangle as the middle open connected component of $\D \setminus ((ab)\cup(bc)\cup(ca))$ (in other words, the connected component of this set that has $a$, $b$ and $c$ on its boundary). A triangle is thus identified with the unordered set of its three  apexes $\{a,b,c\}$.  The set of all hyperbolic triangles will be denoted by $\T$. 

We will denote by $\T_\circ$ the set of all {\em marked} hyperbolic triangles (that can be viewed as the set of ordered triplets $(a,b,c)$ of distinct boundary points that are ordered anti-clockwise on $\partial  \mathbb{D}$). Each hyperbolic triangle corresponds to three marked triangles (one just has to distinguish one apex in order to mark the triangle).

Let us stress the fact that this is a slight abuse of terminology, as our hyperbolic triangles always have their apexes on the boundary of $\D$ (these  triangles  are called \emph{ideal} in hyperbolic geometry, but since we will not use any other triangles in the present paper, we will simply omit to specify that we always mean ideal triangles). Notice also that with our definition, any triangle is open and has non-empty interior. \medskip 

Clearly, it is possible to identify the set $\T_\circ$ of all marked triangles with the group $\M$ of all M\"obius  transformations (hyperbolic isometries) of the unit disk, i.e., the group of transformations  of the type   \begin{eqnarray*} \phi_{z_0, \theta} : z \longmapsto \frac{\mathrm{e}^{ \mathrm{i} \theta} (z-z_{0}) }{\overline{z_0} z -1}, \quad  \mbox{where } z_0 \in \D \mbox{ and }\theta \in [0, 2 \pi). \end{eqnarray*}
Indeed, for each $(a,b,c) \in \T_\circ$, there exists a unique $\phi \in \M$ (that we can therefore call $\phi_{a,b,c}$), such that $\phi ((1,  \mathrm{j},  \mathrm{j}^2)) = (a,b,c)$ (where $ \mathrm{j} =\exp ( 2 \i \pi / 3)$ denotes the cubic root of the unity). Furthermore, in this identification $ \phi_{a,b,c} \leftrightarrow (a,b,c)$, the left-multiplication by an element $\phi$ of $\M$ corresponds  to taking the image of $(a,b,c)$ under this map i.e., $\phi \circ \phi_{(a,b,c)} 
\leftrightarrow (\phi (a), \phi(b), \phi(c))$. 

\medbreak

Let us now describe natural measures that one can define on these different sets.
Recall first that the hyperbolic metric  on $ \mathbb{D}$ is defined by \begin{eqnarray*} \mu &=&  \frac{4 \ \mathrm{d}x  \mathrm{d}y}{ \pi (1-x^2-y^2)^2}  \end{eqnarray*} which is (up to a multiplicative constant) the unique measure on $\D$ that is invariant under the group $\M$. Note that all triangles are equivalent up to hyperbolic isometry, so that they all have the same hyperbolic area. It is easy to check that this area is finite, and we normalized  $\mu$ in such a way that the common area of all triangles is equal to $1$.

The identification of $\T_\circ$ with the locally compact Lie group $\M$ immediately shows  that, up to a multiplicative constant, there exists a unique Haar measure on $\T_\circ$ that is invariant under the group $\M$ (i.e. corresponding to the measures on $\M$ invariant under left-multiplication). In other words, there exists a unique M\"obius-invariant measure $\nu_{\circ}$ on $\T_\circ$ (up to a multiplicative constant). Recall that $\M$ is unimodular, so that $\nu_\circ$ is also invariant under right-multiplication. Here are a couple of simple explicit constructions of $\nu_\circ$:

\begin {itemize}
\item Consider the product measure $\mu \otimes \lambda$ on $\D \times [0, 2 \pi)$, where $\mu$ is the hyperbolic measure in $\D$ and $\lambda$ is the uniform probability measure on $[0, 2\pi)$. Each pair $(z_0, \theta)$ in this set defines the isometry $\phi_{z_0, \theta}$ in $\M$, and it is easy to check that the image measure of $\mu \otimes \lambda$ in $\M$ is invariant under right-multiplication. In other words, one can view $\nu_\circ$ as the image of $\mu \otimes \lambda$ under the map
$ (z_0, \theta) \mapsto \phi_{z_0, \theta} ((1, \mathrm{j}, \mathrm{j}^2))$ (note that with this construction, the point $ \mathrm{e}^{ \mathrm{i}\theta}z_0$ is the ``hyperbolic'' center of the triangle); it is easy to check that indeed  
$$ \mu \otimes \lambda ( \{ (z_0, \theta)  \ : \ 0 \in  \phi_{z_0, \theta} ((1, \mathrm {j}, \mathrm {j^2})) \} )= 1 . $$
\item Another way to construct the M\"obius-invariant measure on $\T_\circ$ goes as follows: Define on $\R^3$ the measure 
$$ \frac {  \mathrm{d}u \  \mathrm{d}v \ \mathrm{d}w }{ \  |u-v| \ |v-w| \ |w-u |},$$ 
where we restrict ourselves to the triplets $(u,v,w)$ that are ordered anti-clockwise around $\partial \H$ (i.e. $u<v<w$, $v<w<u$ or $w<u<v$).  
Clearly this measure is invariant under the transformations $z \mapsto - 1/z$, $z \mapsto \lambda z$ and $z \mapsto z + z_0$ for $z_0 \in \R$ and $\lambda >0$. 
Hence, the image of this measure under $\psi^{-1}$ is a measure on $(\partial \D)^3$ (or rather on $\T_\circ$) that is M\"obius-invariant.  It is therefore necessarily equal to a multiple of $\nu_\circ$. In fact an explicit computation shows that the multiplicative constant is $\pi^2$.
\end {itemize}
  
Similarly, if a measure $\eta$ on $\mathcal T$ is M\"obius-invariant, we can note that the measure on marked triangles obtained by marking one corner uniformly at random among the three, is an invariant measure on ${\mathcal T}_{\circ}$, and therefore a multiple of $\nu_\circ$. It follows that $\eta$ is a multiple of the measure $\nu$ obtained from $\nu_\circ$ by the natural projection from $\T_\circ$ onto $\T$. 

Recall that the $\nu_\circ$-mass (and therefore the $\nu$-mass also) of the set of all triangles that have the origin in their interior is equal to one. By M\"obius invariance, the same holds for the set of all triangles that have a given point $z$ in their interior. In the sequel, $P_z$ (respectively $P_z^\circ$) will denote the probability measure on $\T$ (resp.\,$\T_\circ$) that is equal to $\nu$ (resp.\,$\nu_\circ$) restricted to those triangles that contain $z$. The probability measure $P_0$ will be used repeatedly in the sequel. 
 
\subsection {M\"obius-invariant triangulations}

\label{section:cit}

A  (hyperbolic) triangulation $ \mathbf{T}$ of $\D$ is a disjoint collection of  hyperbolic triangles of $\D$. Since every triangle has non-empty interior, such a collection is finite or countable.
If $\TT$ is a triangulation and $z \in \D$, we define $T(z)$ to be the (unmarked) triangle of $\TT$ that contains $z$ if it exists, and $T(z) = \emptyset$ otherwise. Clearly, if we choose a fixed countable dense family $(z_q)_{q \in Q}$ in $\D$, then the family $(T(z_q), q \in Q)$ fully describes $\TT$.  This gives a way to define a natural sigma-field on the set of all triangulations of $\D$ { which we will implicitly use from now on. Note that this sigma-field in fact does not depend on the choice of the dense family $(z_{q})$ -- indeed if we order $Q$ and identify $Q$ with  $\N$, then for all $z \in \D$, one  has 
$$T(z) = \cup_{j \ge 0} T(z_j) 1_{z \in T( z_j) \mbox{ \ \footnotesize and } z \notin T(z_1) \cup \ldots \cup T(z_{j-1}) }$$ (and where, here and throughout the paper, $T(z_j)1_A $ is equal to the triangle $T(z_j)$ when the event $A$ holds and to the empty set otherwise). It is also easy to check that this sigma-field coincides with that associated with the Hausdorff topology on $ \mathbb{D}$
(but we will not use this fact).}
 
We say that a triangulation $\TT=  (T_j)$ is \emph{complete} if the hyperbolic area of $\D \setminus \cup_j T_j$ is zero. Most triangulations that we will consider in this paper will be complete.
Let us make two side-remarks here (that will not be useful in this paper so that we just mention them leaving the details to the interested reader):
\begin{remark} \label{dense-complete}
 Note that a complete triangulation $\TT$ is dense, in the sense that the union of the triangles of $ \TT$ are dense in $ \mathbb{D}$. However there exist dense triangulations that are not complete (an analogy that one can keep in mind is that there exist open subsets $O$ of $[-1,1]$ that are dense in $[-1,1]$, but with Lebesgue measure that is strictly smaller than $2$ -- such an open subset $O$ then loosely speaking corresponds to the intersection of the interiors of the triangles with $[-1, 1]$). \end{remark}

\begin{remark} A \emph{lamination} is a closed set of $\mathbb{D}$ that can be written as a disjoint union of  hyperbolic lines. One says that a lamination is maximal if it is maximal for inclusion among laminations. It is not hard to see that the complement of a maximal lamination is composed of disjoint open (ideal) triangles and thus is a hyperbolic triangulation, see \cite{Bon98} for more details on hyperbolic laminations.
\end{remark}

We say that a 
 random complete triangulation $\TT$ is {\em M\"obius-invariant} if it is invariant (in law) under the action of each M\"obius transformation of the unit disk. In words, it is M\"obius-invariant if, for any conformal map $\phi$ from $\D$ onto $\D$, $\phi ( \TT)$ and $\TT$ have the same law.
 \medskip 

Suppose now that $\TT$ is such a M\"obius-invariant random complete triangulation. \medskip

By $ \mathcal{M}$-invariance, the quantity $P(T(z) \ne \emptyset)$ is independent of $z \in \mathbb{D}$ and must be equal to  { $1$} by completeness.   \label{loitrig}We can also associate an infinite ``counting'' measure $\eta$ on $ \mathcal{T}$ with  $ \mathbf{T}$ as follows: For any measurable set $A$ of triangles in $ \mathcal{T}$, we define $\eta (A)$ to be the expected value of the number of triangles $T_j \in \TT$ that fall in $A$. By $ \mathcal{M}$-invariance of $ \mathbf{T}$ it follows that $\eta$ is a M\"obius-invariant measure on $\mathcal T$. Note also 
that $\eta ( 0 \in T )$ is equal to the mean number of triangles of $\TT$ that contain the origin, which is equal to $1$ since $\TT$ is almost surely complete. 
Consequently one has $\eta=\nu$, and that every $z \in \mathbb{D}$, $T(z)$ is distributed according to $P_z$.

\medbreak

Of course, it is worth checking if non-trivial M\"obius-invariant triangulations exist at all. 
Here is a construction of the simplest one of all, based on the standard Farey-Ford tiling of $\D$.
Suppose that $ \tau$ is a given (unmarked) hyperbolic triangle. We construct deterministically a triangulation $\mathrm{Ref}(\tau)$ containing $\tau$ by reflections: It is the only triangulation $ \boldsymbol{\mathfrak{T}}$ with the property that for any triangle $T \in \boldsymbol{\mathfrak{T}}$, if  $\phi$ denotes any one of the three M\"obius transformations that map $\{1, \mathrm{j},  \mathrm{j^2}\}$ onto $T$, then $T$ has exactly three adjacent triangles in $\boldsymbol{\mathfrak{T}}$ that are $\phi(\{1,\omega^1,\omega^2\})$, $\phi(\{\omega^2,\omega^3,\omega^4\})$  and $\phi (\{\omega^4,\omega^5,1\})$ where $\omega = \exp(\i\pi/3)$.
 It is elementary to check that $ \mathrm{Ref}(\tau)$ is well-defined and is a complete hyperbolic triangulation. The triangulation $ \mathbf{F}:=\mathrm{Ref}(\{1, \mathrm{j}, \mathrm{j}^2\})$ is called the Farey-Ford tiling, see \cite[Chapter 8]{Bon09}. 
One can identify the set of all marked triangles in $\mathbf {F}$ with the discrete subgroup $G$ of $\M$ that leaves $\mathbf {F}$ invariant. In this way, the set of all triangles of $\mathrm{Ref} (\tau)$ is nothing else than the family of all $g(\tau)$, where $g$ spans $G$ (and because we have been using marked triangles to define $G$, each triangle of 
$\mathrm{Ref} (\tau)$ appears three times in this list).  
  
\begin{figure}[!ht]
 \begin{center}
 \includegraphics[height=4cm]{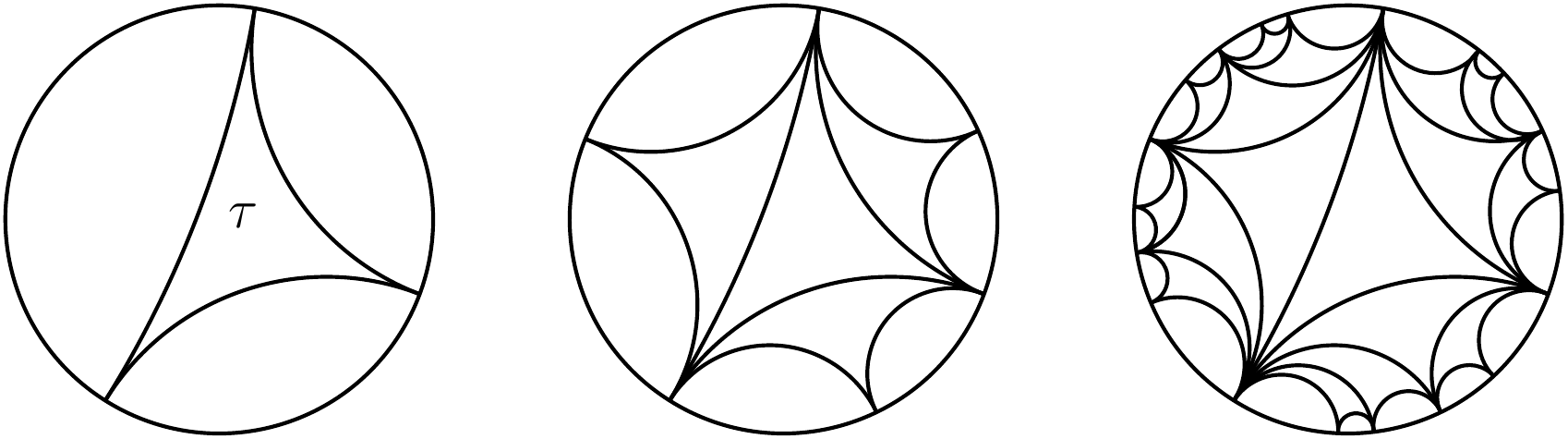}
 \caption{Construction of $\mathrm{Ref}(\tau)$ started from some $\tau$.}
 \end{center}
 \end{figure}
 \begin {proposition} \label{farey-ford}
 If $T_{0}$ is distributed according to $P_0$, then  $\mathrm{Ref}(T_{0})$ is M\"obius-invariant.
\end {proposition} 

\proof 
Note that the knowledge of any triangle in $\mathrm{Ref}(T_{0})$  characterizes the entire triangulation. It therefore suffices to prove that if $T(z)$ is the triangle that contains $z$ in this triangulation, it is distributed according to $P_z$ (as this will imply that the law of the triangulation is invariant under any M\"obius transformation that maps $0$ onto $z$, because $P_z$ is the  image of $P_0$ under such a hyperbolic isometry).

As $\nu$ is invariant under $\M$, it follows that for each $g \in G$, the measures $M_g$ and $M_g'$ defined on the set of pairs of triangles by 
$$M_g ( A) = \nu (\{ T \in  \mathcal{T}\  : \  (T, gT) \in A  \}) \quad \hbox { and } \quad M_g'(A) = \nu ( \{ T'  \in \mathcal{T} \ : \  ( g^{-1} T', T' ) \in A  \} )$$ 
are identical. It follows of course that $\sum_g M_g = \sum_g M_g'$. 
But, if one restricts $\nu$ to those triangles that contain the origin, one obtains $P_0$, and furthermore, almost surely, only one unmarked triangle in $\mathrm{Ref}(T_{0})$ does contain $z$, i.e. $T(z) = gT_0$ for exactly three $g$'s in $G$, and no other $gT_0$'s do contain $z$. Hence, it follows that 
\begin {eqnarray*}
\lefteqn
{3  P ( T(z) \in A )   =   \sum_g \nu ( \{ T \ : \ 0 \in T, z \in gT, gT \in A \})} \\
 && =  \sum_g \nu ( \{ T' \ : \ 0 \in g^{-1} T', z \in T',  T' \in A \} )
  =  3 P_z (A).    
\end {eqnarray*}
\endproof

\subsection {Markovian triangulations}

\label{markov}

Let us now define the additional Markovian property that we will require for our random triangulations.
A first rather weak assumption would be that, conditionally on $T(0)$, the intersection of $\TT$ with the three connected components of $\D \setminus T(0)$ are independent. Mind that the previous randomized Farey-Ford example satisfies this property (indeed, conditionally on $T(0)$, all other triangles are deterministic).

\begin{figure}[!ht]
 \begin{center}
 \includegraphics[height=6cm]{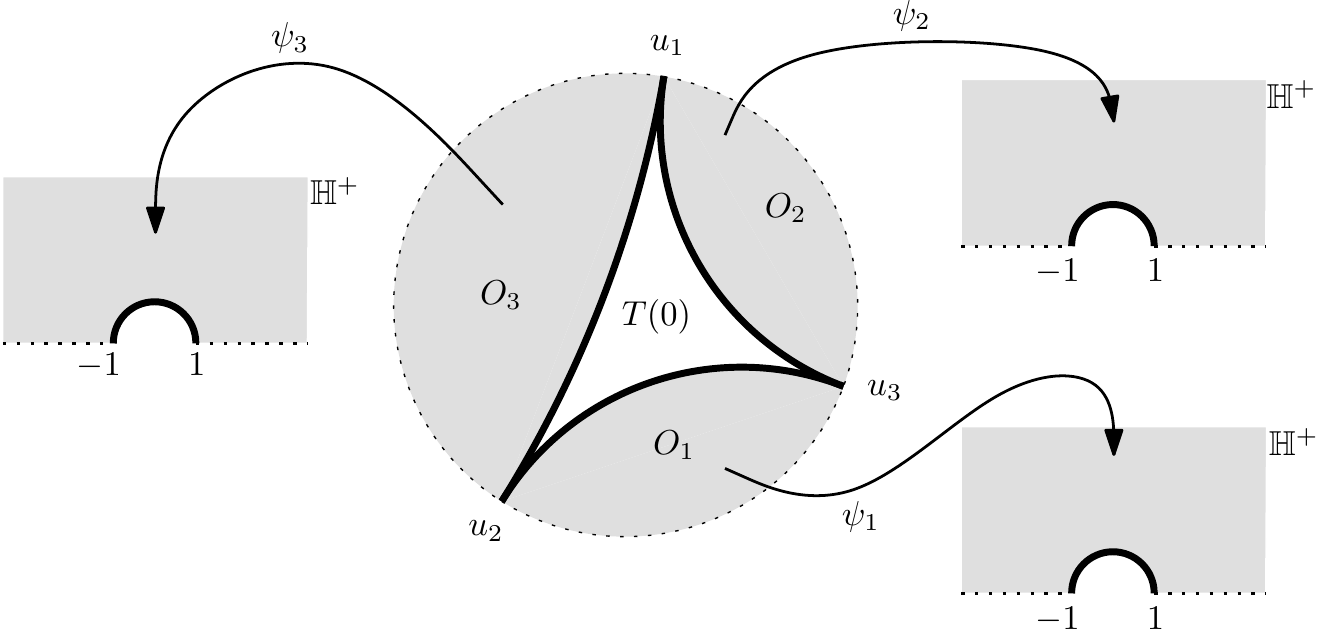}
 \caption{\label{notations}The three maps $\psi_1$, $\psi_2$ and $\psi_3$.}
 \end{center}
 \end{figure}

Our Markovian condition will be stronger. Suppose that we denote the three apexes of $T(0)$ by $u_1$, $u_2$ and $u_3$ ordered anti-clockwise, and the three connected components of $\D \setminus T(0)$ by $O_1$, $O_2$ and $O_3$ in such a way that $u_j \notin \partial O_j$, see Figure \ref{notations}. In order to define which of the three apexes is denoted by $u_1$, we can for instance just choose it at random among the three.

We will say that a random complete triangulation is Markovian if:
$$ 
\hbox {\it Conditionally on $(u_2, u_3)$, $\TT \cap O_1$ is independent of $(\TT\cap O_2, \TT\cap O_3)$.}
$$
Note that (because the triangulation is complete), one can recover $u_1$ from $\TT \cap O_2$ and $\TT \cap O_3$, so that $\TT\cap O_1$ is conditionally independent of $u_1$, given $(u_2, u_3)$.

 Note also that if $\TT$ is also M\"obius-invariant, then the same statement holds for the restriction of $\TT$ to the three connected components of the complement of $T(z)$ (for any given fixed point $z$ in $\D$). 
We are now ready to state our main result:

\begin {theorem} 
\label{main}
There exists exactly one (law of a) Markovian M\"obius-invariant complete triangulation in $\D$. 
\end {theorem}

Until the rest of this section, $\TT$ will denote a random Markovian M\"obius-invariant complete triangulation, and we will start to study its properties. 
Let us  make a first observation. Define for each $j \in \{1, 2, 3\}$ a conformal transformation $\psi_j$ from $\D$ onto $\H$   that maps $O_j$ onto the domain
 $$\H^+ : = \{ z \in \H \ : \ |z| > 1 \},$$ see Fig.\,\ref{notations}.
 We choose a way to define $\psi_j$ that is a deterministic function of $O_j$ (so that $\psi_1$ does not depend on $u_1$ etc.) --  let us for instance pick $\psi_1$ so that $\psi_1 (u_3) = -1$, $\psi_1 (u_2) = 1$ and $|\psi_1' (u_2)| = 1$). 
In this way, each $\widetilde \TT_j := \psi_j ( \TT \cap O_j)$ is a triangulation of $\H^+$.

The following statement is a consequence (i.e., a reformulation)
 of our Markovian assumption for M\"obius-invariant complete triangulations. 

\begin {lemma}
\label {Lsigma} The variables 
$T(0)$, $\widetilde \TT_1$, $\widetilde \TT_2$, $\widetilde \TT_3$ are independent, and the latter three have the same distribution. Furthermore, this common distribution $\sigma$ is invariant under the one-dimenional group of all M\"obius transformations $\phi$ of $\H$ such that $\phi(\H^+)= \H^+$.
\end {lemma}

\proof
Let us (for notational convenience) decide that each triangle of $\TT$ has been marked at random and independently (this therefore defines $u_1(z), u_2(z), u_3(z), O_1(z), \ldots$ for almost all $z \in \D$ because the triangulation is complete in such a way that  if $T(z)=T(z')$ then $u_{i}(z) = u_{i}(z'), O_i(z)=O_{i}(z'), \ldots $ ). We will denote by $\TT_j^z$ the triangulation $\TT$ restricted to $O_j(z)$, $\widetilde{ \mathbf{T}}_{j}^z$ its image in $ \mathbb{H}^+$, and $\theta_1(z)$ will denote the harmonic measure at $z$ of the part of $\partial \D$ between $u_2=u_{2}(z)$ and $u_3=u_{3}(z)$ that does not contain $u_1=u_{1}(z)$. In particular, $2 \pi \theta_1(0) \in (0, \pi)$ is simply the angle at the origin of the triangle $u_2 0 u_3$.

Note that $\widetilde \TT_1, \widetilde \TT_2, \widetilde \TT_3$ are conditionally independent given $T(0)$ (because each $\widetilde \TT_j$ is a deterministic function of $\TT \cap O_j$). In order to derive the full independence, it therefore suffices to check that (for each given $j$), $\widetilde \TT_j$ and $T(0)$ are independent. By symmetry, it is sufficient to consider the case $j=1$. As $\widetilde \TT_1$ and $u_1$ are conditionally independent given $(u_2, u_3)$, it is enough to show that $\widetilde \TT_1$ and $(u_2, u_3)$ are independent. 
Because of rotational invariance (and because $\widetilde \TT_1$ does not change if one rotates $\TT$), it finally suffices to show that $\widetilde{\TT}_1$ and $\theta_1(0)$ are independent.

Let $F$ denote a measurable bounded real-valued function on the set of triangulations, and $h$ a measurable bounded function on $\R$. Then,
\begin {eqnarray*}
 E \big[ F( \widetilde \TT_1 ) h (\theta_1 (0) )\big] 
&=&  E \big[F( \widetilde \TT_1 ) h (\theta_1 (0) ) \times \int_{\D} \mathrm{d} \mu(z) 1_{z \in T(0)}\big] \\
&=&  \int_\D \mathrm{d} \mu (z) E \big[ F (\widetilde \TT_1^z ) h ( \theta_1 (0)) 1_{0 \in T(z)} \big] \\
&=&  \int_\D \mathrm{d} \mu (z) E \big[ F (\widetilde \TT_1^0) h ( \theta_1 (z)) 1_{z \in T(0)}\big] \\
&=&  E \big[ F ( \widetilde \TT_1)\big] \times \int_\D h (\theta_1 (z)) 1_{z \in T(0)} \mathrm{d}\mu(z),
\end {eqnarray*}
where we have used the facts that the $\mu$-area of $T(0)$ is one, that $0 \in T(z)$ if and only if $z \in T(0)$, that the triangulation is M\"obius-invariant (and in particular under the hyperbolic isometry that interchanges $z$ and $0$), and finally that $\int_{T(0)} h (\theta_1 (z))  \mathrm{d}\mu (z)$ is a constant that does not depend on the triangle $T(0)$ (because of M\"obius invariance of all quantities involved). This therefore completes the proof of the fact that $\widetilde{\TT}_1,\widetilde{\TT}_2$ and $\widetilde{\TT}_3$ are independent and independent of $T(0)$. They clearly have the same law that we denote by $\sigma$. 
 
It remains to show that $\sigma$ is invariant under all M\"obius transformations of $\H$ that map $ \mathbb{H}^{+}$ onto itself. 
Recall from Section \ref{loitrig} that we know explicitly the distribution of $T(0)$ which has a smooth density with respect to the Lebesgue measure on $(\partial \D)^3$, and we have just seen that $T(0)$ is independent of $\widetilde \TT_1$. Hence, we can say that for any given triangle $(a,b,c)$ that contains the origin, the conditional distribution of $\widetilde \TT_1$ given $T(0)=(u_1,u_2,u_3) =  (a,b,c)$ is $\sigma$. Suppose that $\Phi: \D \to \D$ is some M\"obius transformation, and define $\TT' = \Phi(\TT)$. If we combine the previous decomposition with the M\"obius invariance of $\TT$, we see that for any  $a$, $b$ and $c$, the conditional distribution of $\widetilde \TT_1'$ given $T'(0)= (\Phi(a), \Phi(b), \Phi(c))$ is still $\sigma$ as long as this new triangle contains $0$. 

Let $\phi$ be a fixed M\"obius transformation of $\H$ onto itself such that $\phi(1) = 1$ and $\phi (-1) = -1$. By choosing $a$, $b$ and $c$ appropriately (for instance $a=-1$, $b=-\i \exp (\i \varepsilon))$ and $c= \bar b$ with $\varepsilon$ very small so that the latter two points are very close to $-\i$ and $\i$), we can make sure that if we define 
$\Phi = \psi_1^{-1} \circ \phi \circ \psi_1$ (where $\psi_1$ is the isometry defined deterministically using $b,c$ only, that maps these two points onto $1$ and $-1$, chosen with the same rule as the one we used to define $\tilde \TT_1$) then both triangles 
$(a,b,c)$ and $(\Phi(a), \Phi(b), \Phi(c))$ contain the origin. Furthermore, we see that for this particular triple, when $T(0) = (a,b,c)$, the conditional law of 
$\widetilde \TT_1'$ is that of $\phi (\widetilde \TT_1)$. It follows that the law of $\widetilde \TT_1$ is indeed invariant under $\phi$. 
\endproof

{

We can note that this proves in particular that we can define (in terms of $\sigma$) the conditional distribution of $T(z)$ given $T(0)$. 
We can also note that by M\"obius invariance of $\TT$, for each given $z$, the previous lemma also yields (using the conformal map that swaps $0$ and $z$) a description of the conditional law of the three triangulations corresponding to $\TT$ restricted to each of the three connected components of the complement of $T(z)$, in terms of $\sigma$.  

\medbreak

The next lemma shows that in order to prove uniqueness (in law) of M\"obius-invariant complete Markovian triangulations, it suffices to prove that all their two-dimensional marginals are uniquely determined:}

\begin {lemma}
\label {twodimmar}
If for each $z$, $z' \in \mathbb{D}$  we know the joint law of $(T(z), T(z'))$, then we know the law of the entire triangulation $\TT$.
\end {lemma}

\proof { 
The law of $\TT =( T(z), z \in \D )$ is characterized by the law of its finite-dimensional marginals i.e.\,by the law of $T(Z):=\{T(z_{1}), ... , T(z_{n})\}$ for all finite sets of points $Z= \{ z_1, \ldots, z_n  \}$ in $\D$ with rational coordinates (see the discussion on sigma-fields at the beginning of Section \ref{section:cit}). 

We say that a finite collection of disjoint triangles in the unit disk is good if each connected component of the complement of the union of these triangles in the disk has at most two neighboring triangles in this collection. The left picture of Figure \ref {fidi} represents a set $T(Z)$ that is not good, because the shaded component neighbors three different triangles of $T(Z)$. 

However, because $\TT$ is almost surely complete, for each given $Z$, with probability one,  it is possible to add to $T(Z)$ finitely many triangles of $\TT$ in order to turn it into a good collection. In fact, there is a minimal way to do this, and we call $\tilde T (Z)$ the corresponding finite collection of triangles of $\TT$. Note that (for each given $Z$) $T(Z)$ is a deterministic function of $\tilde T(Z)$ (it consists of those triangles in $\tilde T(Z)$ that contain a point of $Z$), so that the distribution of $\tilde T (Z)$ contains all the information about the distribution of $T(Z)$.  

When $Z'$ is a finite set of points in $\D$, we say that a finite collection $ \mathfrak{T}$ of triangles is in ${\mathcal A} (Z')$ if $\mathfrak{T}$ is good, and if each triangle in $\mathfrak{T}$ corresponds exactly one point of $Z'$ (i.e. each triangle of $\mathfrak{T}$ contains exactly one point of $Z'$ and each point of $Z'$ is in a triangle of $\mathfrak{T}$).  In particular, if $\TT$ is our random triangulation, the event $\{T(Z') \in \mathcal A (Z')\}$ holds iff $T(Z')$ is good and if each triangle of $T(Z')$ contains exactly one point of $Z'$. 

Note that if $Z$ is some other given finite family of points, by looking at $T(Z')$ only, we can see whether $\tilde T (Z) = T(Z')$. Similarly, we can also check whether this event holds or not by looking at $\tilde T(Z)$ only (it suffices to check that $\tilde T(Z) \in \mathcal {A} (Z')$). 
Suppose that for a given $Z'$, we know the law of $T(Z') 1_{T' \in {\mathcal A} (Z')}$. Then, clearly, for each given $Z$, we know the law of 
$$T(Z') 1_{T(Z') \in {\mathcal A} (Z') \hbox  { and } T(Z')= \tilde T (Z) }= \tilde T(Z) 1_{\tilde T (Z) \in \mathcal {A} (Z')}.$$

But if $I$ denotes the family of finite sets $Z'$ with rational coordinates, then for any given $Z$, 
$$P ( \tilde T(Z) \in \cup_{Z' \in I} {\mathcal A} (Z')) = 1$$ 
(because each triangle of $\tilde T(Z)$ contains some point with rational coordinates). 
Hence, it follows that if, for all $Z'$, one knows  the law of  $T(Z') 1_{T(Z') \in {\mathcal A} (Z')}$, then one can reconstruct the law of $\tilde T(Z)$ and therefore that of $T(Z)$.

Finally, it remains to prove that for each finite $Z$ (we use $Z$ instead of $Z'$ now), the law of $T(Z)1_{T(Z) \in {\mathcal A} (Z)}$ is fully determined by the knowledge of all two-dimensional marginal distributions of $(T(z), T(z'))$ for $z,z' \in \mathbb{D}$. We are going to do one further reduction step:
Let us now suppose that for some finite $Z$ we have $T(Z) \in {\mathcal A} (Z)$ (recall that 
this means that $T(Z)$ is good, and that each triangle of $T(Z)$ corresponds exactly to one point of $Z$). 
This defines naturally a connected tree structure $G$ on $Z$, where each $z_j$ has one, two or three neighbors in the graph, see Figure \ref {fidi}. We can therefore decompose the 
event $\{T \in {\mathcal A} (Z)\}$ according to the tree-structure that $T(Z)$ induces on $Z$. Hence, it suffices to describe for each $Z$ and each possible tree-structure $\Gamma$ on $Z$, the law of $T(Z) 1_{T \in {\mathcal A} (Z, \Gamma) }$, where 
$$ {\mathcal A}(Z, \Gamma ) = {\mathcal A} (Z) \cap \{ G = \Gamma \}.$$

\begin{figure}[!ht]
 \begin{center}
 \includegraphics[height=5cm,angle=0]{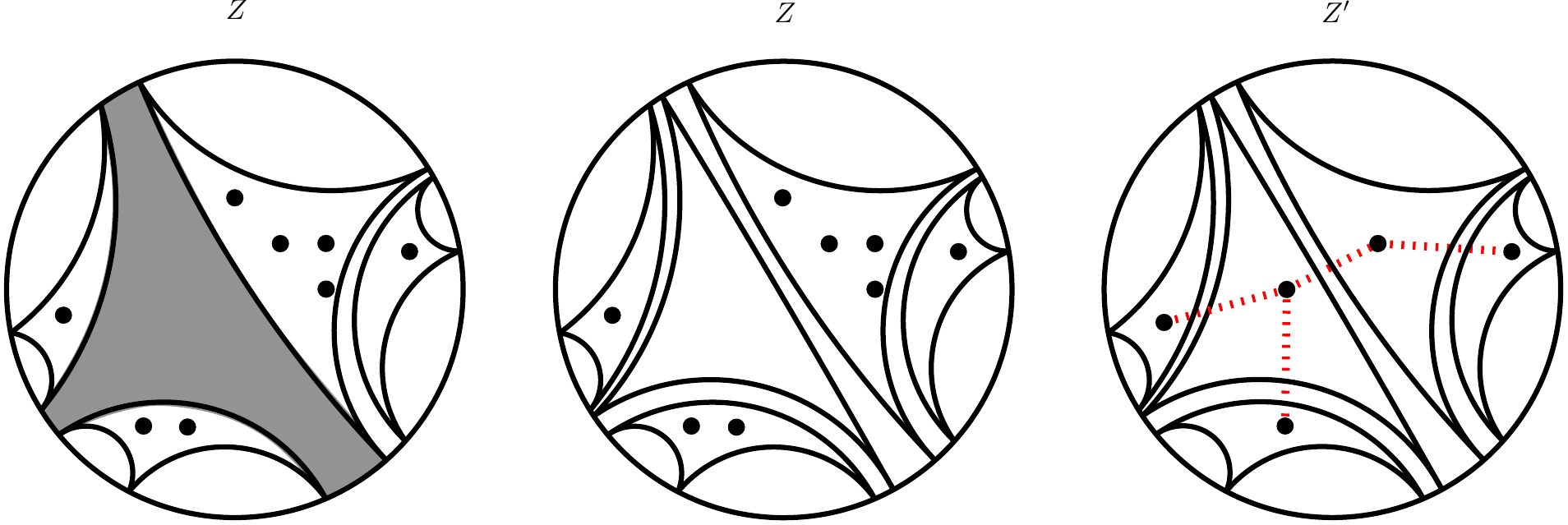}
 \caption{\label{fidi} A configuration $T(Z)$, its completed configuration $ \tilde{T}(Z)$ and a configuration in $ \mathcal{A}(Z')$.}
 \end{center}
 \end{figure}

We are going to proceed by induction on the number of points in $Z$.
Suppose now that we know the law of all two-dimensional marginals $(T(z), T(z'))$, and that for each $Z$ with no more than $n$ points, and for each tree-structure $\Gamma$ on $Z$, we know the law of $T(Z) 1_{T \in {\mathcal A} (Z, \Gamma) }$. Let us show that, we then know it also for all $Z= \{z_1, \ldots, z_{n+1} \}$ with $n+1$ points and all tree-structure $\Gamma$ on $Z$. 
Let us choose such a $Z$ with $n+1$ points and a tree-structure $\Gamma$ on $Z$. Consider a leaf-point (i.e., a point in $Z$ with just one $\Gamma$-neighbor) -- by relabeling the points,  we can assume that this leaf is $z_1$ and that its $\Gamma$-neighbor is $z_2$. Our assumptions and previous results show that we know:
\begin {itemize}
\item The distribution of $(T(z_2), \ldots,  T(z_{n+1}))$, when restricted to the event that it defines the tree-structure obtained by removing the leaf $z_1$ from $\Gamma$.
\item The fact that conditionally on $T(z_2)$, on the event that it separates $z_1$ from the other $n-1$ points, $T(z_1)$ is independent from $(T(z_3), \ldots, T(z_n))$.
\item The joint distribution of $(T(z_1), T(z_2))$ (and therefore also the conditional distribution of $T(z_1)$ given $T(z_2)$).
\end {itemize}
This shows readily that we know the distribution of $T(Z) 1_{T \in {\mathcal A} (Z, \Gamma) }$: Indeed, first sample $(T(z_2), \ldots, T(z_{n+1}))$, look if it is compatible with ${\mathcal A} (Z, \Gamma)$, and then sample $T(z_1)$ according to the conditional distribution given $T(z_2)$. 

Hence, we have proved our claim by induction over $n$, which  provides a characterization of the law of all $T(Z) 1_{T \in {\mathcal A} (Z, \Gamma) }$, and therefore by our previous arguments, of $\TT$ itself.
}\endproof    

\section{Uniqueness}
\label {S3}

\subsection {Warm-up}

In order to help those readers who are not so acquainted with the theory of regenerative sets, we briefly review some very classical facts on this topic (we refer to \cite {Ber96,Ber99} for details). Those readers who are familiar with these objects can safely skip this subsection. \medskip 

Suppose that we are given a random non-empty closed subset $F$ of $\R_+$  such that almost surely, $0 \in F$, $F$ is not bounded, and the Lebesgue measure of $F$ is  equal to $0$. Suppose furthermore that it satisfies the following regenerative property: For any $t \ge 0$, if we define $X_t = \min [t, \infty) \cap F$, then the law of $F_t:= (F \cap [X_t, \infty)) - X_t $ is equal to that of $F$. We also assume that $F_t$ is independent of $(X_t, F \cap [0, t])$. 
In the standard terminology, this means that $F$ is a ``light'' (because its Lebesgue measure is $0$) regenerative subset of $\R_+$. 

Then, for each given small positive $x$, we can discover the intervals of length greater than $x$ in $\R_+ \setminus F$ from left to right. This defines (at least for small enough $x$) a sequence $\xii^x:= (\xii_1^x, \xii_2^x, \ldots )$ of lengths. The previous assumptions readily imply that this a sequence of independent identically distributed random variables, that have some common law $\rho_x$. 
Furthermore, when $x' < x$, the fact that $\xii^{x'}$ is almost surely a subsequence of $\xii^x$ implies that 
$\rho_x = \rho_{x'}( \cdot  \mid [x,\infty) )$.
Hence, we can define a measure $\rho$ on all of $(0, \infty)$ with the property that for all small enough $x$, 
$$ \rho_x ( A ) = \frac {\rho (A \cap [x, \infty))} {\rho ([x, \infty))}.$$
The measure $\rho$ is unique up to a multiplicative constant and is in a way describing the relative likelihood of appearance of intervals of a certain length in the complement of $F$. Note that it can happen that the total mass of $\rho$ is infinite, which corresponds to the fact that there can be infinitely many (small) intervals in the complement of $F \cap [0,1]$ say. 

Now, it turns out that the measure $\rho$ completely characterizes the law of the random set $F$. For instance, we can define simultaneously for each $x >0$, a sample of $\xii^x$ in such a way that they are all compatible (i.e. $\xii^{x}$ is almost surely a subsequence of $\xii^{x'}$ when $x' < x$). Then, the left-hand point of the interval corresponding to $\xii_n^x$ will be the sum of all intervals (of arbitrary length) that have appeared before it, which can be recovered from the knowledge of all $\xii^{x'}$ for $x' < x$.   
 
One convenient way to express this is to use a Poisson point processes: This is a random countable collection $\mathcal {P}:= (t_i, x_i)_{i \in I}$ in $\R_+ \times \R_+$ where we introduced an artificial time-parametrization at which the intervals appear.  Intuitively, the existence of the point $(t_i, x_i)$ in ${\mathcal {P}}$ means that at time $t_i$, an interval of length $x_i$ appears. If for a given $x$, we write down the sequence of lengths $x_i$ greater than $x$ in their order (with respect to time) of appearance, one gets a sample of $\xii^x$.  Then, the position of the left-point of the interval corresponding to $i_0$ can be recovered from $\mathcal {P}$ as it is equal to
$ \sum_{ i \in I \ : \ t_i < t_{i_0}} x_i$.

Note that if we were  looking at the set $\tilde F:=  \{ \exp(t) \ : \ t \in F \}$ instead of $F$, we would have had a set satisfying similar properties to $F$: With obvious notation, the law of $\tilde F_t / \tilde X_t$ is equal in law to $\tilde F_0$, and that in order to recover $\tilde F$ from $\rho$, one has to replace the sum of all jumps $x_i$ by the multiplication of all $\exp (x_i)$. 
We shall use rather natural generalizations of these ideas in the next subsection. 

\subsection{Accordion and Poisson point process}

\label{uniqueness}
\begin{figure}[!ht]
 \begin{center} 
 \includegraphics[height=5cm]{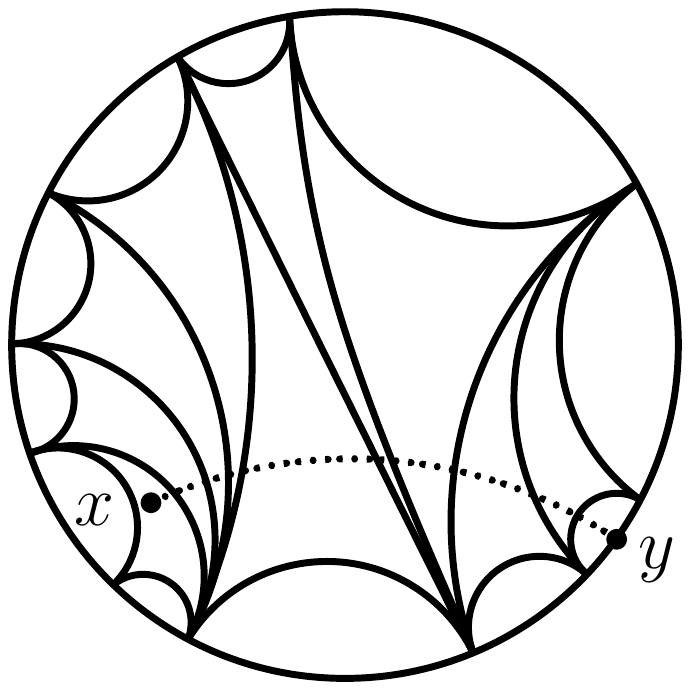}
\caption{ \label{accordion} An accordion.}
 \end{center}
 \end{figure}

 Let $\TT$ be a  complete triangulation of $\D$.  For $x \not= y$ in $\overline{\D}$, we define the \emph{accordion} between $x$ and $y$ in $\mathbf{T}$ as the collection of all triangles $T\in \mathbf{T}$ intersecting the part of the hyperbolic line between $x$ and $y$ that goes through these two points,
  and denote it by $ \mathrm{Acc}_{\mathbf{T}}(x,y)$, see Fig.\,\ref{accordion}.

Suppose now that $\TT$ is  M\"obius-invariant complete and Markovian.
Clearly, if we know the distribution of the accordion $ \mathrm{Acc}_{ \mathbf{T}}(0,1)$, this will characterize the law of $(T(0), T(u))$ for all $u \in (0,1)$ and therefore (by M\"obius invariance and  Lemma \ref {twodimmar}) also the distribution of $\TT$. The goal of the present section will be to show that there is (at most) one possible law for this accordion. The basic idea will be to see that it necessarily corresponds to some particular subordinator (when one ``discovers'' the triangles of the accordion from $0$ towards $1$). In the next section, we shall check that the random triangulation defined using these random accordions is indeed M\"obius-invariant complete and Markovian.

For notational convenience, we will now choose to work in the upper half-plane $\H$ instead of the unit disk $\D$. For the remainder of this section, $\TT$ will denote a random M\"obius-invariant complete Markovian triangulation of the upper half-plane (i.e. the image under $\psi$ of such a triangulation of $\D$).
Note that almost surely, $\infty$ and $0$ are on the boundary of none of the triangles of $\TT$ (this follows from rotational invariance and the fact that the set of triangles is countable -- one can also just look at the second characterization of the measure $\nu$ in the preliminaries). In other words, all triangles are bounded and bounded away from the origin.
We are going to focus on the accordion between $ \mathrm{i}$ and $\infty$ in $ \mathbf{T}$, that will be denoted by $\mathrm{Acc}_{\TT} ( \mathrm{i}, \infty)$ (we will omit to specify that we are working in $\H$).

\begin{figure}[!ht]
 \begin{center}
 \includegraphics[height=6cm]{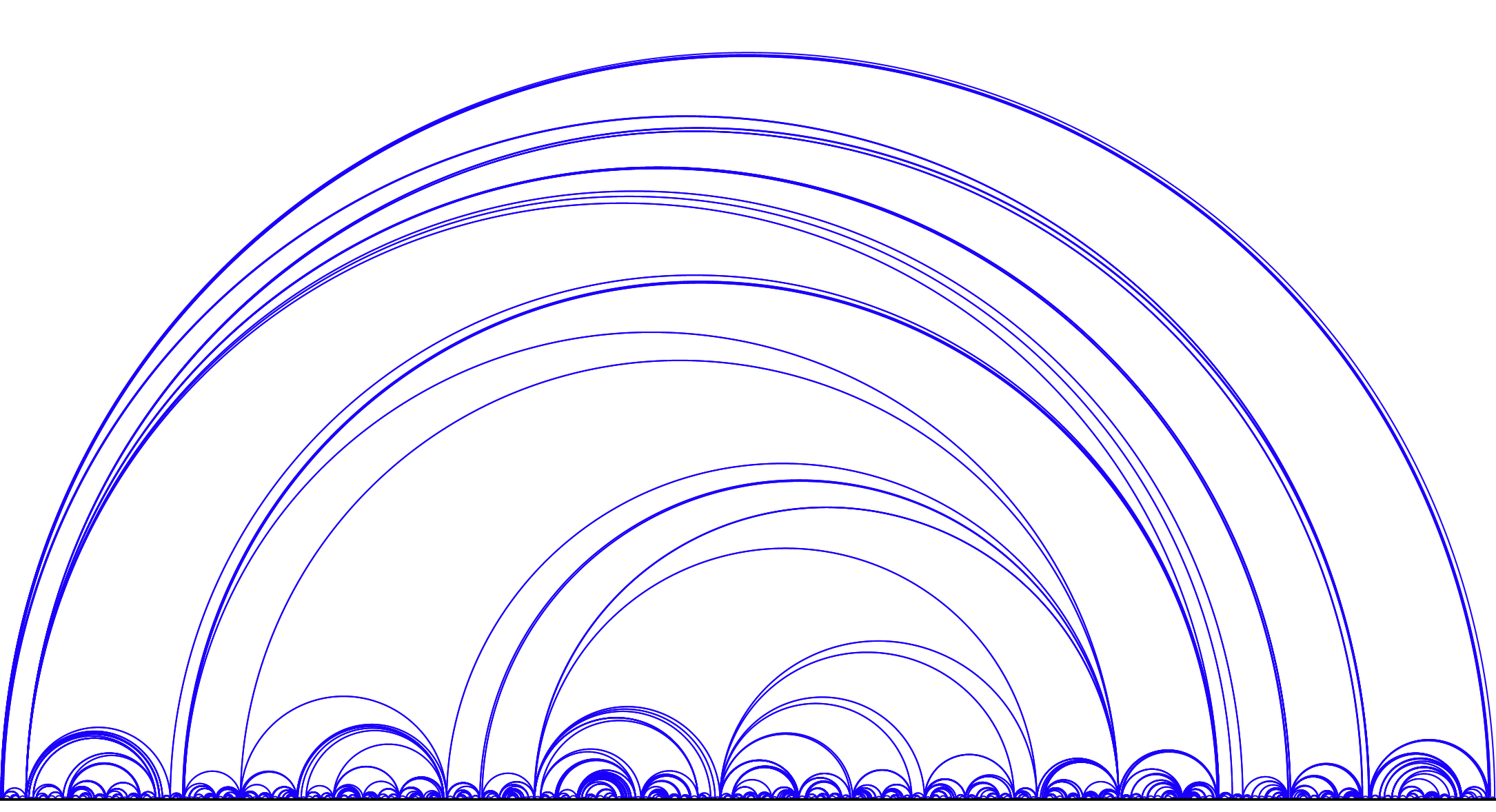} \hspace{1cm}
 \caption {Sample of (part of) our triangulation in $\H$.}
 \end{center}
 \end{figure}
 
Almost surely, for all { positive rational} $y$, one of the boundaries of $T( \mathrm{i}y)$ does separate $ \mathrm{i}y$ from $\infty$. We denote it by $(\ell_y r_y)$ where $\ell_y < 0 < r_y$. 
Clearly, $-\ell_y$ and $r_y$ are non-decreasing functions of $y$. We can therefore define $(\ell_y, r_y)$ for  all positive $y$ simultaneously (including those that are in no triangle) by choosing the right-continuous version of $y \mapsto ( \ell_y, r_y)$.

{
Let us first outline the idea of the proof: If we discover the triangle $T( \mathrm{i} y)$, the conditional law of the part of the triangulation that is ``above'' this triangle can be described via $\sigma$, and it is (modulo taking its image under the affine map that maps $\ell_1, r_1$ onto $-1, 1$) always the same. It follows from this observation that the closure of the set 
$ \{ (r_y - \ell_y ) / (r_1 - \ell_1) \ : \ y > 1 \}$ is the exponential of a regenerative set, just as $\tilde F$ in the end of the warm-up subsection. It can therefore be described thanks to a Poisson point process -- the fact that the triangulation is complete will imply that the Lebesgue measure of this set $\tilde F$ is $0$. The set $\tilde F$ does however not contain enough information in order to reconstruct the accordion because when a triangle appears, one needs to know which one of the two processes $\ell$ or $r$ is jumping. We will therefore describe the accordion via a slightly enriched Poisson point process that contains this additional information. \bigskip

For all positive $y$, define
 \begin{eqnarray*}\widetilde y 
 &:=& \sup\{ v \in (0, y)  \   : \  T( \mathrm{i}v) \not= T( \mathrm{i}y) \}. \end{eqnarray*}
For all $y$ such that $T( \mathrm{i}y) \not= \emptyset$, the third vertex of $T(\i y)$ 
(apart from $\ell_y$ and $r_y$) is necessarily one of the two points $\ell_{\widetilde y-}$ or $r_{\widetilde y-}$. Note also that the jumps of the process $(\ell,r)$ exactly correspond to the triangles of $\op{Acc}_{\mathbf{T}}( \mathrm{i},\infty)$
 (i.e., the set $\J$ of ``jumping heights'' is equal to $\{ \widetilde y \ : \ y > 0  \}$). Note also that $y$ can never be a simultaneous jumping height for $\ell$ and $r$ (because almost surely, no $T( \mathrm {i} y)$ is a quadrilateral).

For each positive $y$, we denote by $\varphi_y$ the affine map that maps $(\ell_{y-}, r_{y-})$ onto $(-1, 1)$. We can then describe the jumps of $T(iy)$ by defining for each $y \in \J$,
$X(y)$ to be the image of the third apex  of $T(\i y)$  (apart from $\ell_{y-}$ and $r_{y-}$) under $\varphi_{y}$. 
In other words, 
 \begin{eqnarray} \label{jump} X(y) &=& \varepsilon(y) \times \left(2\left( \frac {r_y - \ell_y}{r_{y-} - \ell_{y-}} \right)-1 \right)  \end{eqnarray}
where
 \begin{eqnarray*}\varepsilon(y) &=& 1_{\{r_y \not= r_y-\}} - 1_{\{ \ell_{y} \not= \ell_{y-} \} }. \end{eqnarray*}
 When $y \notin \J$, we can declare $X(y)$ to be equal to an abstract cemetery point $\delta$.

{ Note that for all  $y_1 >1$, and all $x > 1$, the number of jumps $\tilde{y}$'s in $[1, y_1]$ such that $|X( \tilde y)| > x$ is finite. Hence, it follows that the { collection $(X( \tilde{y}) 1_{|X(\tilde{y})| > x}, \tilde{y} \ge 1)$} almost surely defines an ordered discrete sequence  $\zeta^x = (\zeta_1^x, \zeta_2^x, \ldots)$ in $\R \setminus [-x,x]$. Note that when $x' > x$, the sequence $\zeta^{x'}$ is a deterministic subsequence of $\zeta^x$. We define $\CC$ to be this nested family of sequences $(\zeta^x, x >1 )$ (we can not view it just as one sequence, because infinitely many ``small'' jumps occur before any given jump).  

An equivalent way to encode $\CC$ is to define it as the process of jumps $(X(\tilde{y}), \tilde{y} \geq 1)$ but defined modulo increasing time-reparametrization i.e., only the order of arrivals of the jumps matters. 

In the following, for $x \geq1$ we let $I_{x} = \mathbb{R} \backslash [-x,x]$. In particular $I_1 = \mathbb{R}\backslash[-1,1]$.
}

\begin {lemma}
\label {32} The ordered (but unparametrized) set of jumps   $ \mathcal{C}$ has the same distribution as the ordered family of  jumps (modulo increasing time-reparametrization) of a Poisson point process $ \mathcal{P} = \{(t_{i},x_{i})_{i \in I}\}$ on $\R_{+} \times \R$ with intensity $ \mathrm{d}t \otimes \rho$, where $\rho$ is some sigma-finite measure on $\R \setminus [-1,1]$.
\end {lemma}
\proof {  Let $t >1$. By the spatial Markov property applied to the triangle $T( \i t)$ (i.e., the push-forward of Lemma \ref {Lsigma} by $\psi$), one deduces that the ordered (but unparametrized) collection of jumps $(X({ \tilde{y}}), { \tilde{y}} \geq t)$ is independent of $(X({ \tilde{y}}), 1 \leq { \tilde{y}} \leq t)$ and has the same distribution as the ordered family of jumps  $(X({ \tilde{y}}), { \tilde{y}} \geq 1)$. 

Fix $x_{0} >1$ such that there almost surely exists a jump in $I_{x_{0}}$. We deduce from the above remark that for every $1<x<x_{0}$,  the discrete random sequence $(X(\tilde{y})\mathbf{1}_{|X(\tilde y)|  \in I_{x}}, \tilde{y}>1)$ has the same distribution as i.i.d.\,samples from a certain probability measure $\rho_{x}$ on $I_{x}$. Furthermore, the distributions $\rho_{x}$ satisfy the compatibility 
 \begin{eqnarray*} \rho_{x'}\big(  . \,\left|\, I_{x}\big)\right. &=& \rho_{x}  \end{eqnarray*}
for all $1 < x' < x$. Consequently, on can define uniquely a sigma-finite measure $\rho$ on $ I_{1}=\mathbb{R} \backslash [-1,1]$ such that $\rho( . \cap I_{x}) / \rho(I_{x}) = \rho_{x}$ and $\rho(I_{x_{0}}) =1$. It is then easy to see that the jumps of $ \mathcal{C}$ have the same distribution as the unparametrized jumps of a Poisson point process $\mathcal{P} = \{(t_{i},x_{i})_{i \in I}\}$ on $\R_{+} \times \R_{}$ with intensity $ \mathrm{d}t \otimes \rho$, see the warm-up section. Details are left to the reader.

The sigma-finite measure $\rho$ has an arbitrary multiplicative normalization (but note that the multiplicative constant does not change the law of the ordered family of jumps, it just changes the time-parametrization). }
\endproof 

We have now seen that if a M\"obius-invariant complete Markovian triangulation $\TT$ exists, then one can associate with it a measure $\rho$ that describes the law of the jumps of $\mathcal C$, and we have also seen that the distribution of $T(\i)$ is necessarily the image of $P_0$ under $\psi$. 
Furthermore, the Markovian property (Lemma \ref{Lsigma}) shows that $T(0)$ and the jumps of $\mathcal C$ are independent.
The following lemma proves that conversely one can  recover the law of $\mathrm{Acc}_\TT (\i,\infty)$ from $\rho$ and $T(i)$:
\begin {lemma}
\label {31}
The distributions of $T( \mathrm{i})$ and $\CC$ fully characterize the law of  $\mathrm{Acc}_\TT ( \mathrm{i}, \infty)$. 
\end {lemma}
\proof It is clear that $ \mathrm{Acc}_{ \mathbf{T}}( \i, \infty)$ can be recovered from the two processes $(\ell_{y})_{y \geq 1}$ and $(r_{y})_{y \geq 1}$ and the initial triangle $ {T}( \i)$.  { More precisely, instead of the full processes $(\ell)$ and $(r)$,  it suffices to know $(\ell, r)$ \emph{up to time reparametrization} to reconstruct the accordion. Indeed, only the range of $(\ell, r)$ matters in order to define  $ \mathrm{Acc}_{ \mathbf{T}}(\i,\infty).$}

We first claim that the ranges of the processes $(\ell)$ and $(r)$ are both of zero one-dimensional Lebesgue measure. Recall that the triangulation $\TT$ is almost surely complete, and M\"obius-invariant, so that any given point in $\H$ is almost surely in the interior of some triangle of $\TT$. Hence, the (one-dimensional) Lebesgue measure of the intersection $\mathcal I$ of the imaginary line with the closure of the union of all arches $(\ell_y r_y)$ is almost surely equal to $0$. 
 
Indeed, assume that the Lebesgue measure of the intersection of the range of $(\ell)$ with some interval $[-l_1, -l_2]$ is positive. Clearly, one can associate to each point $\ell_y$ of this range, a point $\i y$ on the imaginary half-line, in such a way that for any $y<y'$, $|y' - y| > |\ell_{y'} - \ell_y |/ K(l_1)$ for some constant $K(l_1)$. Hence, it follows readily  that the one-dimensional Lebesgue measure of $\mathcal I$ is positive. Since this is prohibited, we conclude that  the range of $(\ell)$ (and of $(r)$, by the same argument) is 
almost surely  of zero (one-dimensional) Lebesgue measure.

{ As $(\ell)$ and $(r)$ are monotone functions, the fact that their ranges are of zero Lebesgue measure implies that the range of $(\ell, r)$ is characterized by its jumps (which themselves are described by $ \mathcal{C}$) and by its initial value (which is given by $T(\i)$). Hence, we can recover, up to time reparametrization, the process $(\ell, r)$ from $ \mathcal{C}$ and $T(\i)$. This is sufficient to reconstruct $ \mathrm{Acc}_{ \mathbf{T}}(\i,\infty)$.}\endproof

\subsection{Identification of the jump measure}

It now remains to show that (up to a multiplicative constant), there is in fact at most one possibility for the measure $\rho$ defined in Lemma~\ref{32}. This will follow from the M\"obius invariance of the measure $\sigma$  as  heuristically described in the introduction. 

Let us suppose for the remainder of this section that $\TT$ is a M\"obius-invariant complete Markovian triangulation, and that $\rho$ and $\mathcal {P}$ are defined as in Lemma \ref {32} (and that $\mathcal {P}$ is coupled with $\mathcal C$ in such a way that their ordered family of jumps are identical).  { If $\varphi$ is a M\"obius transformation of $ \mathbb{H}$, the action of $\varphi$ can be extended to boundary $ \partial \mathbb{H} = \mathbb{R}$. We will implicitly use this extension is what follows. Note that it is sufficient for an arch or a triangle to track down its apexes $\in \partial \mathbb{H}$ to know it entirely.}
\begin {lemma}
\label {33}
The image measure of $\rho$ under any M\"obius transformation of the upper half-plane into itself that fixes $-1$ and $1$ is proportional to $\rho$.
\end {lemma}

{ \proof  Fix $x_{0}>1$ in such a way that $\rho ( \{x_0, -x_0 \} ) = 0$ and $\rho(I_{x_{0}})>0$. Define $\ell_1< 0 < r_1$ as before so that $(\ell_1, r_1)$ is the top boundary of $T (\i)$. We denote $\tilde{x}$ the first jump of the point process $ \mathcal{P}$ such that $|\tilde{x}|>x_{0}$ and write $\tilde{\ell} < 0 < \tilde{r}$ for the feet of the bottom hyperbolic line of the triangle corresponding to the jump $\tilde{x}$.  Let $\widetilde a$ denote the 
third apex of this triangle. See Figure\,\ref{setup}.

\begin{figure}[!ht]
 \begin{center}
 \includegraphics[height=6cm]{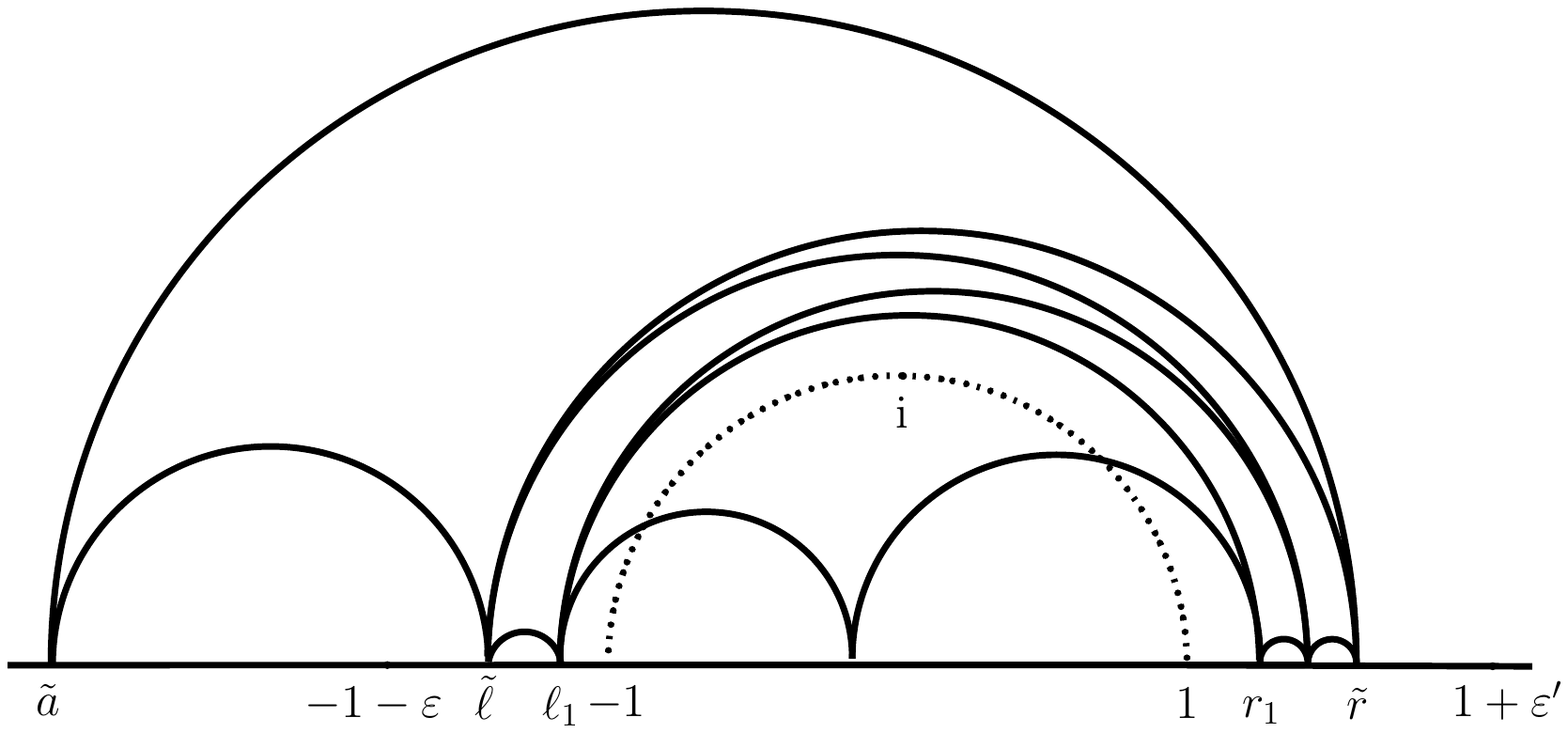}
 \caption{ \label{setup}Setup of the proof.}
 \end{center}
 \end{figure}
For each small $\varepsilon$ and $\varepsilon'$, we consider the events
$$ A( \varepsilon, \varepsilon') :=  \{ -1-\varepsilon< \widetilde \ell <\ell_1 < -1 \hbox { and } 1 < r_1 <  \widetilde r < 1 + \varepsilon' \}.$$

  By standard properties of Poisson point processes the event $A(\varepsilon,\varepsilon')$ is independent of $\tilde{x}$ which is distributed according the measure $\rho_{x_0}$. Thus for any Borel positive $f :  \mathbb{R}\backslash [-1,1] \to \mathbb{R}_{+}$ with compact support  we have \begin{eqnarray} \label{conv1} E\left[ f( \tilde{x}) \mid A(\varepsilon,\varepsilon')\right] &= & \rho(I_{x_0})^{-1} \int_{I_{x_{0}}} \rho (\mathrm{d}a) f(a)  \end{eqnarray} We will now let $x_{0} \to 1$. To avoid subsequent normalizations, we consider another positive measurable function $g : \mathbb{R}\backslash [-1,1] \to \mathbb{R}_{+}$ with compact support: Using the last display and letting $x_0 \to 1$ we have
  \begin{eqnarray} \label{conv2} \frac{E\left[ f(\tilde{x}) \mid A(\varepsilon,\varepsilon')\right]}{E\left[ g(\tilde{x}) \mid A(\varepsilon,\varepsilon')\right]} & \xrightarrow[ \varepsilon, \varepsilon'\to0]{}  &  \frac{\int_{ I_1} \rho (\mathrm{d}x) f(x)}{\int_{ I_1} \rho (\mathrm{d}x) g(x)}  \end{eqnarray}

  On the other hand, $ \tilde{x}$ can be related to the geometric quantity $ \tilde{a}$ as follows. When $\varepsilon$ and $\varepsilon'$ are both very small (and $A(\varepsilon, \varepsilon')$ holds) then  the jump $\tilde{x}$ is necessarily very close to   the first foot $ \tilde{a}$ of the accordion with absolute value larger than  $x_{0}$. Thanks to the above remark, the $\tilde{x}$ can be replaced by the geometric $ \tilde{a}$ in the left-hand side of \eqref{conv2}.

Let us now suppose that $\varphi$ is a M\"obius transformation that maps $\H$ onto itself with $\varphi (-1) = -1$ and $\varphi (1) = 1$. Note in particular that since the semi-circle $(-1, 1)$ is preserved by $\varphi$, for every $\varepsilon,\varepsilon'>0$ there exist $\delta,\delta'>0$ such that if $A(\delta,\delta')$ is satisfied for $\TT$ then $A(\varepsilon,\varepsilon')$ holds for $\varphi(\TT)$ and furthermore $T (\i) = T( \varphi( \i))$). Since $\varphi ( \TT)$ and $\TT$ are identically distributed it follows readily using the same arguments as before that 

\begin{eqnarray} \label{conv4} \frac{E\left[ f(\varphi(\tilde{a})) \mid A(\varepsilon,\varepsilon')\right]}{E\left[ g(\varphi(\tilde{a})) \mid A(\varepsilon,\varepsilon')\right]} & \xrightarrow[ \varepsilon, \varepsilon'\to 0]{} &  \frac{\int_{ I_1} \rho (\mathrm{d}x) f(x)}{\int_{ I_1} \rho (\mathrm{d}x) g(x)}  \end{eqnarray} Thus comparing (\ref {conv4}) with \eqref{conv2} (with $\tilde x$ replaced by $\tilde a$) we deduce that
  \begin{eqnarray*} \frac{\int_{ I_1} \rho (\mathrm{d}x) f(\varphi(x))}{\int_{ I_1} \rho (\mathrm{d}x) g(\varphi(x))} &=& \frac{\int_{I_1} \rho (\mathrm{d}x) f(x)}{\int_{I_1} \rho (\mathrm{d}x) g(x)}.  \end{eqnarray*}
Hence the push-forward of $\rho$ under the map $\varphi$ is indeed a multiple of $\rho$.}
}
\endproof

A natural candidate for the measure $\rho$ is the measure $\zeta$ on $\R \setminus [-1, 1]$ defined by 
$$ \zeta (\mathrm{d}x) = \frac {2 \mathrm{d}x }{|x|^2 - 1} 1_{\{|x| > 1 \}}$$
as it is the only measure (up to a multiplicative constant) that is invariant under all M\"obius transformations of $\H$ that fix the boundary points $-1$ and $1$
(note for instance that it is the image of the measure  $ \mathrm{d}x / x$ on $\R_+$ under the map $(1+x)/(1-x)$).
Indeed:

\begin {lemma}
\label {34}
The measure $\rho$ defined in Lemma \ref {32} is necessarily equal to a constant times the measure $\zeta$.  
\end {lemma}

\proof
Let us consider $\varphi(x) = (1+x) / (1-x)$ as above.  Clearly, if $\widetilde \rho$ denotes the push-forward of $\rho$ under $\varphi^{-1}$, 
this measure $\widetilde \rho$ on $\R_+$ will satisfy the property that the image of $\widetilde \rho$ under any map $z \mapsto \lambda z$ for positive $\lambda$ (these are the M\"obius transformations of $\H$ onto itself that fix $0$ and $\infty$) is a multiple (that may depend on $\lambda$) of $\widetilde \rho$. It follows that for some real $\alpha$ and some positive constant $c$, 
$$ \widetilde \rho (\mathrm{d}x) = c x^{-\alpha} 1_{\{x > 0 \}} \mathrm{d}x.$$
We want to show that $\alpha$ is necessarily equal to $1$. Let us assume that $\alpha < 1$. Then, $\tilde \rho [\varepsilon, \infty) = \infty$ while $\tilde \rho [0, \varepsilon) < \infty$ for any $\varepsilon>0$. In terms of $\rho$, this implies in particular that 
$ \rho [1, \infty) <  \infty$. But the proof of Lemma \ref {31} then tells us that the set of $(r_y)_{y \ge 1}$ has no accumulation points i.e. that all $r_y$'s are isolated. In particular, this implies that if $E$ denotes the set of all corners of triangles in $\TT$ that separate $0$ from $\infty$ in $\H$, then $E_+ :=  E \cap \R_+$ is almost surely discrete in the sense that for all $0< a < b $, $E_+ \cap [a,b]$ is finite (here $[a,b]$ denotes the horizontal segment between $a$ and $b$). 

On the other hand, for any $\varepsilon>0$, $\rho [-1-\varepsilon, -1] = \infty$. This readily shows there are infinitely many jumps for $(\ell_y)$ while $(r_y)$ only jumps finitely many times. In particular, we see that almost surely, there exist $b<a<0$ such that there are infinitely many points in the intersection of $E_- := E \cap \R_-$ with the horizontal segment $[b,a]$.     

Finally, because of invariance of the law of $\TT$ under the transformation $z \mapsto - 1/z$, we note that $E_+$ has the same law as $\{ -1/z, z \in E_- \}$, which contradicts the previous facts that we just proved for $E_+$ and $E_-$. 
 
We therefore conclude that $\alpha \ge 1$. In exactly the same way, we can exclude the possibility that $\alpha > 1$ (because then $\rho [1, \infty) = \infty$ while $\rho (-\infty, -1] < \infty$). Hence, we see that $\rho$ is a multiple of the image under $\varphi$ of $x^{-1}\mathrm{d}x 1_{\{x >0\} }$, i.e., a multiple of $\zeta$.
\endproof

The previous lemmas therefore describe the joint law of $(T(\i), T( \i y))$ for any given $y >1$. 
But, for any $z$ and $z'$ in $\D$, there exists some $y \ge 1$ and a M\"obius transformation from $\D$ onto $\H$ that maps $z$ onto $\i$ and $z'$ onto $ \i y$; by M\"obius invariance, we can therefore describe the joint law of $(T(z), T(z'))$, and by Lemma \ref {twodimmar}, we have completed the proof of the uniqueness part of Theorem \ref{main}:
\begin {proposition}
There exists at most one (law of a) complete M\"obius-invariant Markovian triangulation.
\end {proposition}

\section {Existence}
\label {S4}

The goal of this section is to define the candidate for the random triangulation, and to check that it is complete, Markovian and M\"obius-invariant.

\subsection {The half-plane accordion} \label{construction}

In order to define a random accordion in $\H^+$ (i.e., what will turn out to be our distribution $\sigma$), we start with a Poisson point process
 $\mathcal{P}= \{(t_{i},x_{i})_{i \in I} \}$  on $\R_{+} \times (\R \backslash [-1,1])$ with intensity $\mathrm{d}t \otimes \zeta $.

{ We then construct two pure jump processes $(L_{t})_{t \geq 0}$ (for left) and $(R_{t})_{t \geq 0}$ (for right) that jump only on the jumping times of $\mathcal{P}$ whose jumps (defined as in \eqref{jump}) are the $x_{i}$'s. 
 Set $L_{0} = -1$ and $R_{0}=1$. The idea is that $L$ is decreasing, that $R$ is increasing, and that when a  jump  $(t,x)$ occurs, then $R-L$ is multiplied by $(|x|+1)/{2}$, and that $L$ jumps only if $x<-1$ and $R$ jumps only if $x >1$. 

More precisely, if we set $\P_t := \{ ( t_i , x_i) \in \mathcal{P} \ : \ t_i \le t \}$, then we can first define 
 \begin{eqnarray*}\Delta_t &: =&  (R_0 - L_0) \prod_{(t_i, x_i) \in \P_t} \frac{|x_i|  +1}{2}\end{eqnarray*}
which is the exponential of the pure jump process with intensity given by the image of $\zeta$ under the mapping $x \mapsto \log \frac{|x|+1}{2}$. }
It is easily checked that this subordinator is well-defined (that it does not blow up), using the explicit expression for its jump measure.
Then, we simply set  
\begin{eqnarray}
R_{t} &:=& R_0 + \sum_{ (t_i, x_i) \in \P_t } (\Delta_{t_i} -  \Delta_{t_i-}) 1_{x_i > 1} \\ 
L_{t} &:=& L_0 - \sum_{ (t_i, x_i) \in \P_t } (\Delta_{t_i} -  \Delta_{t_i-}) 1_{x_i < -1} . 
\end{eqnarray}
As $\Delta_t =  R_{t}-L_{t}$ is almost surely finite for all $t$,  the two processes $(L_{t})_{t\geq 0}$ and $(R_{t})_{t\geq 0}$ are well-defined as well (note that $L$ is non-increasing and that $R$ is non-decreasing). Note also that $R_t \to \infty$ and $L_t \to -\infty$ almost surely as $t \to \infty$. 

{ 
Equivalently, for each $i \in I$, we can write
$$ 
(L_{t_i}, R_{t_i}) = \varphi_{t_i-}^{-1} ((-1, x_{t_i})) 1_{ x_{t_i} > 1} +  \varphi_{t_i-}^{-1} ((x_{t_i}, 1)) 1_{x_{t_i} < -1} $$
where $\varphi_{t}^{-1}$ denotes the affine map that maps $(-1, 1)$ onto $(L_{t-}, R_{t-})$.
}

We are now ready to define our accordion in $\H^+$. For each $(t_i, x_i) \in \mathcal {P}$, define the hyperbolic triangle in $\H$ with three corners given by 
$(L_{t_i}, R_{t_i-}, R_{t_i})$ if $x_i > 1$, and by $(L_{t_i}, L_{t_i-}, R_{t_i})$ if $x_i < -1$. The definition clearly ensures that each of these triangles is separating the semi-circle $\{ z \in \H \ : \ |z| = 1 \}$ from $\infty$ in $\H^+$ and that these triangles are disjoint.

\begin{figure}[!ht]
\begin{center}
\includegraphics[width=12cm]{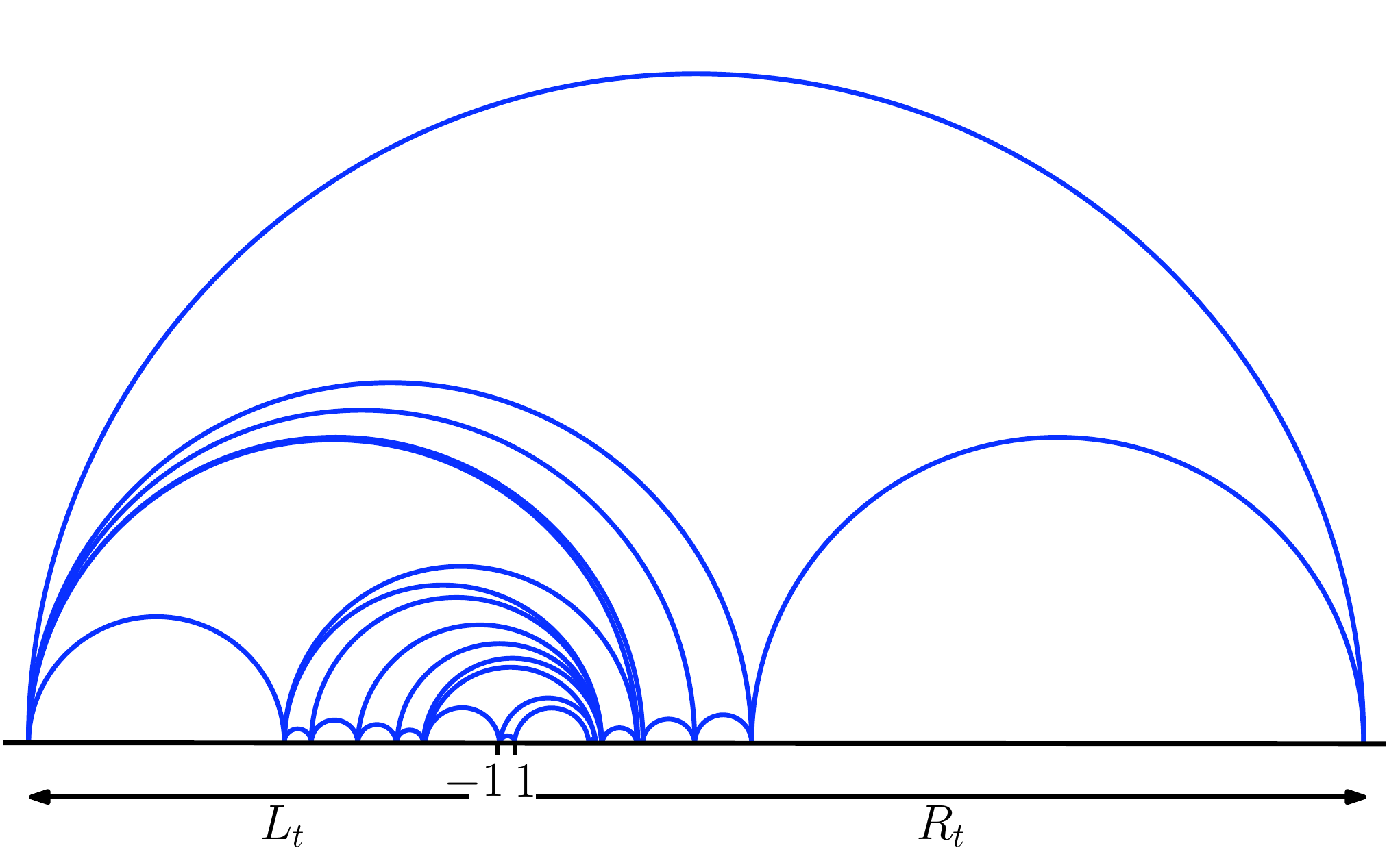}\\ 
\caption{ \label{figaccord} Sample of a piece of $\mathcal{A}(\mathcal{P})$. }
\end{center}
 \end{figure}
 
In fact, in order to indicate the fact that this accordion is from the semi-circle $(-1,1)$ to $\infty$ in $\H$, we will denote it by $\A (\P)_{[(-1, 1) \to \infty; \H ]}$ (and omit the $\H$ when it is clear that we are working in $\H$, and then simply write $\A (\P)_{(-1,1) \to \infty}$).

\label{existence} 
We can immediately extend our construction to the case where the initial position $(L_0, R_0) = (l_0, r_0)$ (for $r_0 > l_0$) 
is different than $(-1, 1)$ and this defines the accordion $\A (\P)_{(l_0,r_0) \to \infty}$. 
It is easy to check that that this new accordion has the same distribution as the the image of 
$\A (\P)_{(-1, 1) \to \infty}$ under the linear map that maps $-1$ onto $l_0$ and $1$ onto $r_0$. 

In other words, we have in fact defined $(L_t, R_t)$ as a Markov process on $\{ (l,r) : l < 0 < r\}$  with translation-invariant and scale-invariant transition kernel (the process started from $(-1+x, 1+x)$ has the same law as $(x+L_t, x+R_t)_{t \ge 0}$ when $L_0 = -1$ and $R_0 = 1$, and on the other hand, the process started from $(-r, r)$ has the same law as $(rL_t, rR_t)_{t \ge 0}$).   

\subsection {Towards completeness}

Let us now prove the following statement: 
\begin{lemma} \label{box}
Almost surely, the ranges of $(L_{t})_{t\geq 0}$ and  $(R_{t})_{t\geq 0}$ restricted to any compact interval of $\mathbb{R}$ are both of box-counting dimension $0$.
\end{lemma}
\proof Consider an auxiliary subordinator $(Q_t, t \ge 0)$ defined by
$$ Q_{t} =\sum_{ (t_{i},x_{i})  \in \P_t } \frac{x_i-1}{2}   \mathbf{1}_{x_{i}>1}.$$
 Let $T>0$. It is clear from its construction that the process $(R_t, t\ge 0)$ jumps exactly when  the process $(Q_t, t \ge 0)$ jumps, and that up to time $T$, { the size of a jump 
 of $R$ is less than the corresponding jump of $Q$ multiplied by $\Delta_T$,  indeed if $t<T$ is a jump time for $R$ we have $$R_{t}-R_{t-} = \Delta_{t}-\Delta_{t-} = \Delta_{t-}\frac{x_{i}-1}{2} = \Delta_{t-}(Q_{t}-Q_{t-}).$$}
 Thus,  the box-counting dimension of the range $R[0,T]$ is almost surely 
 not larger than that of $Q [0,T]$ (because the former set is the image of the latter under a Lipschitz map). But the box-counting dimension of $Q[0,T]$ is easily seen to be  almost surely equal to zero  (see \cite[Chapter 5.1.1]{Ber99} or \cite[Chapter III.5]{Ber96}, and use the behavior of $\zeta( \mathrm{d}x)$ near $x=1$). As the process $L$ has the same law as $-R$, the lemma follows.
\endproof

Similarly as in Lemma \ref{31}, we will translate the previous result on the range of $L$ and $R$ into a property on the set of points of the accordion that are on the imaginary axis. More precisely, let us define the set $ \mathcal{I}_\infty$ of points of the type $\i y$ for $y \ge 1$ that are not  inside  a triangle of $\A (\P)_{(-1,1) \to \infty}$. Then:
\begin {corollary}
\label{Kinfini}
 The (one-dimensional) Lebesgue measure of $ \mathcal{I}_\infty$ is almost surely equal to zero. 
\end {corollary}
\proof For any $k >1$, let us define the set $J_k$ of points $\i y$ for $y \ge 1$ that are in the closure of the union of the semi-circles $((L_t, R_t), t <\sigma_k)$, where $\sigma_k$ is the first time at which $\max (R_t, -L_t) \ge k$. 
Clearly, it is sufficient to prove that for any given $k$, this set $J_k$ has almost  surely zero Lebesgue.  
For each positive $\varepsilon$, define  $N_{\varepsilon}$ to be the minimal number of intervals of length $\varepsilon$ that are needed to cover 
$R[0,\sigma_k]\cup L[0, \sigma_k]$.  Lemma \ref{box} in particular implies that almost surely, $N_\varepsilon / \varepsilon^{-1/3}$ vanishes as $\varepsilon$ tends to $0$. 

\begin{figure}[!ht]
  \begin{center}
  \includegraphics[height=5cm]{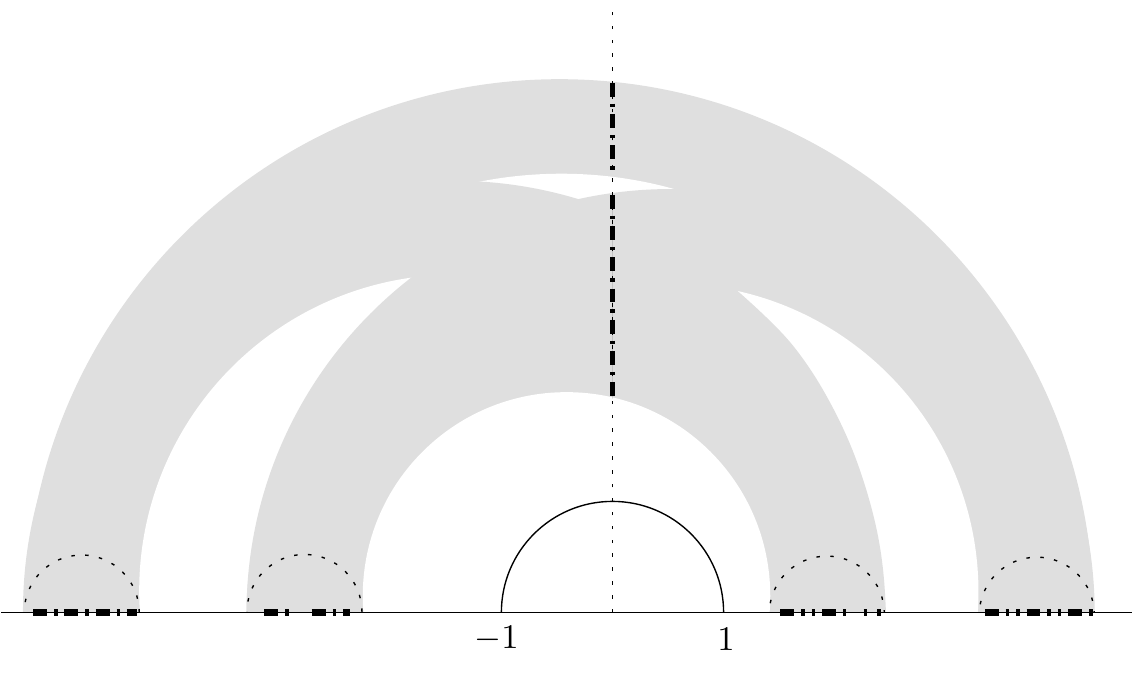}
  \caption{ \label{fig:recouvrement} Sketch of the covering.}
  \end{center}
  \end{figure}

Suppose that $\i y$ (for $y  > 1$) is in no triangle of $\A (\P)_{(-1,1) \to \infty}$. Then it means that one can find one of the intervals of length $\varepsilon$ covering the range of $L$ (let us call it $I_l$), and one of the intervals of length $\varepsilon$ covering the range of $R$ (that we call $I_r$) such that $\i y$ is in one of the semi-circles joining a point in $I_l$ to a point in $I_r$. See Fig.\,\ref{fig:recouvrement}. But, for a given $k$, and any two such intervals, the length of the set of points on the imaginary axis that can be reached in this way is bounded by a constant $C=C(k)$ times $\varepsilon$.  
Thus we have that
$$
 \mathrm{Leb}_{1}(J_k) \leq C \varepsilon  \times N_{\varepsilon}^2 .
 $$
The right-hand side goes almost surely to $0$ as $\varepsilon$ vanishes which concludes the proof of the corollary.
 \endproof

\subsection {Target-independence}
 
Let us now recall a simple classical lemma (see for instance \cite{Ber96} Section O.5.- it can be viewed as a direct consequence of the ``compensation formula'') that roughly states that if we start with a Poisson point process, and modify it in a way that preserves both the independence and the intensity measure then the law of the modified point process is still the same:  
\begin{lemma}[Modification of Poisson point processes]  \label{modif} Let $\mathcal{P}= \{(t_{i},x_{i}), i \in I\}$ be a Poisson point process 
 on $ \mathbb{R}_{+} \times \R$ of intensity $ \mathrm{d}t \otimes  \rho$ (where $\rho$ denotes some measure on $\R$). Let $(H_{t})_{t\geq 0}$ be a predictable process taking values in the space of nonnegative measurable functions  $\R \to \R$, such that almost surely, for every $t\geq 0$ the push-forward of the measure $\rho$ by $H_{t}$ is $\rho$. Then $\mathcal{P}' := \left\{(t_{i},H_{t_{i}}({x_{{i}}})), i \in I \right\}$ has the same law as $\mathcal{P}$.
 \end {lemma}
 
We will use this lemma in order to derive a target-independence property for our accordion  $\A (\P)_{(-1,1) \to \infty}$.  \medskip { 

Fix $a  \in \mathbb{R} \backslash [-1,1]$ and let us define $\psi_{a,0}$ to be the M\"obius map from $\H$ onto itself that maps $(-1, 1, \infty)$ onto $(-1, 1, a)$. We {\em define} the accordion $\A (\P)_{(-1,1) \to a }$ in $\H$ to be the image of $\A ( \P)_{(-1, 1) \to \infty}$ under $\psi_{a,0}$. We finally denote by $\left.\A ( \P)_{(-1, 1) \to \infty\wedge a}\right.$ the sub-accordion of $\A ( \P)_{(-1, 1) \to \infty}$ whose triangles intersect the line $(0a)$ and similarly we denote $\left.\A ( \P)_{(-1, 1) \to a\wedge \infty}\right.$ the sub-accordion of $\A ( \P)_{(-1, 1) \to a}$ whose triangles intersect the line $(0\infty)$.

\begin{proposition}[Target independence] \label{target}For any $a \in \mathbb{R}\backslash [-1,1]$, the two accordions $\left.\A ( \P)_{(-1, 1) \to a\wedge \infty}\right.$ and $\left.\A ( \P)_{(-1, 1) \to \infty\wedge a}\right.$ have the same law.
\end{proposition}
\proof The idea is to decompose the global action of the composition of an accordion with a M\"obius transformation into an iteration of infinitesimal transformations of the jumps by (predictable) functions. Consider $(L)$ and $(R)$ the two functions associated with a standard standard accordion $[(-1,1)\to \infty ; \mathbb{H}]$ and introduce the disconnection time of  $a$ and $\infty$:
 \begin{eqnarray*}\theta_a &:=& \inf \{t \ge 0 \ : \ L_t < a < R_t \}. \end{eqnarray*}

For each jumping time $t_{i} \le \theta_a$, the jump $x_{i}$ in the accordion $[(-1,1)\to \infty]$ corresponds to the image of $L_{t_{i}}$ or $R_{t_{i}}$ under the affine map that sends $(L_{t_{i}-},R_{t_{i}-},\infty)$ onto $(-1,1,\infty)$. Let us now see what is the corresponding jump $x_i^{(a)}$ in the image of $\A ( \P)_{(-1, 1) \to  \infty}$ by $\psi_{a,0}$ that we consider as an accordion growing towards $\infty$, at least as long as $t_i \leq \theta_a$. The jump $x_i^{(a)}$ corresponding to $x_i$ in $\A ( \P)_{(-1, 1) \to a}$  (which is $\psi_{a,0} (\A ( \P)_{(-1, 1) \to  \infty})$ by definition) is the image of $\psi_{a,0}(L_{t_{i}})$ or $\psi_{a,0}(R_{t_{i}})$ under the hyperbolic isometry that sends $(\psi_{a,0}(L_{t_{i}-}),\psi_{a,0}(R_{t_{i}-}),\infty)$ onto $(-1,1,\infty)$. We deduce that $x_i^{(a)}$ is the image of $x_{i}$ by the hyperbolic isometry
 \begin{eqnarray*}\psi_{a,t_{i}} := 
\varphi_{(\psi_{a,0}(L_{t_{i}-}), \psi_{a,0}(R_{t_{i}-} ), \infty) \to (-1,1,\infty)}\circ \psi_{a,0}\circ \varphi_{(-1,1,\infty) \to (L_{t_{i}-},R_{t_{i}-},\infty)}. \end{eqnarray*}


Note that the measure $\zeta$ is invariant under $\psi_{a,t_{i}}$ and that this is a predictable function (with respect to the natural filtration defined by the Poisson point process). When $t > \theta_a$, we simply define $\psi_{a,t}$ to be the identity. Hence, we  deduce from Lemma \ref{modif} that the two ordered but unparametrized families $\{ x_i, i \in I \}$ and $\{  \psi_{a,t_{i}}(x_i), i \in I \}$ have the same law.  This, together with the fact that the jumps characterize the accordion (Lemma \ref{31}), tells us precisely that up to the first time at which one disconnects $\infty$ from $a$, the two accordions 
$\A (\P)_{(-1,1) \to a }$ and $\A ( \P)_{(-1, 1) \to \infty}$ are identically distributed. Note that the final ``jump'' (i.e. the triangle that disconnects $a$ from $\infty$) is also included in this description. \endproof }

 It is therefore possible to couple the two accordions aiming at $a \in \mathbb{R}\backslash[-1,1]$ and $\infty$, in such a way that they coincide up to the  triangle disconnecting  $a$ and $\infty$. This compatibility shows that it is in fact possible to couple accordions (all based on $(-1, 1)$) aiming at all points (with rational coordinates, say) in $\R \setminus [-1, 1]$ in such a way that any two of them coincide up to the first triangle that disconnects their two targets. We define by 
$ \A_{(-1,1), \H}$ the union of all the triangles in this ``accordion tree''. 
Then:
\begin {itemize}
\item
The distribution of $\A_{(-1,1), \H}$ is invariant under the one-dimensional family of  conformal maps from $\H$ onto itself that fix $(-1, 1)$. This follows just from the definition of the accordion targeting other points than $\infty$ via M\"obius invariance.

\item 
This triangulation $\A_{(-1,1), \H}$ is almost surely complete, i.e. almost surely, the two-dimensional Lebesgue measure of the complement in $\H^+$ of the union of all the triangles of $\A_{(-1,1), \H}$ is equal to zero. This is just due to the fact that for any $a \in \R \setminus \{0 \}$, the intersection of this set with the hyperbolic line joining $0$ to $a$ has almost surely zero (one-dimensional) Lebesgue measure (which again follows from the result for $a = \infty$ i.e. from Corollary \ref {Kinfini}, via M\"obius invariance). 
\end {itemize}

This target independence is reminiscent of the ``locality property'' of SLE$_6$ \cite {LSW}.

\subsection{Reversibility of the accordion}

On top of being invariant under  M\"obius transformations of $\H$  that leave $ \mathbb{H}^+$ invariant,  the measure $\zeta$ possesses another property that will yield ``reversibility'' of the accordion. Recall from the end of Section \ref{existence} that we can view $M_t=M_t^+ = (R_t, L_t)_{t \ge 0}$ as a pure jump Markov process on the space $\{ (r,l) \in \R^2 \ : \ l <0<  r \}$. 

We will also use the Markov process $M_t^-= (L_t^-, R_t^-)$ defined on the same state space, but aiming to $0$ instead of to $\infty$. It is defined exactly as $M_t$ except that the measure $\zeta$ is replaced by the measure $\zeta^- (\mathrm{d}x) = 2  \mathrm{d}x / (1-x^2)$ with support in $[-1, 1]$. Note that $\zeta^-$ is the image of $\zeta$ under the map $x \mapsto -1/x$, so that it  follows (using the same arguments as in the proof of Proposition \ref{target}) that if $M_0=M_0^-= (-1, 1)$, then the two processes $(L_t^-, R_t^-)_{ t\ge 0}$ and $ (-1/R_t, -1/L_t)_{t\ge 0}$ have the same law. 

Let us define, for any $u< v$ in $\R$, the two measures $\zeta_{[u,v]}$ and $\zeta_{[v,u]}$ that are supported respectively on $[u,v]$ and $\R \setminus [u, v]$ with respective densities
$$ 
\zeta_{[u,v]} (\mathrm{d}w) = \frac {(v-u) \mathrm{d}w}{ (v-w) (w- u)} 
\hbox { and } 
\zeta_{[v,u]} (\mathrm{d}w ) = \frac {(v-u) \mathrm{d}w}{(w-v) (w-u)}.
$$
Note that these two measures are invariant under any M\"obius transformation of $\H$ that fixes $v$ and $w$, and that they are the push-forward of the measure $ \mathrm{d}x / x$ on $\R_+$ by any M\"obius 
transformation from $\H$ onto itself that maps $0$ and $\infty$ on $u$ and $v$ (in that order for $\zeta_{[u,v]}$ and in the reversed order for $\zeta_{[v,u]}$; in fact, the measure $\mathrm{d}x /x$ on $\R_+$ can be interpreted as $\zeta_{[0, \infty]}$). Hence, all these measures are images of each other under some hyperbolic isometry.
Note also that $\zeta_{[1, -1]}$ is exactly our measure $\zeta$ and that $\zeta_{[-1, 1]}  = \zeta^-$.

The definition of our Markov process $M_t$ shows that $\zeta_{[v,u]}$ describes its jump intensity measure (i.e. the location of the new point after the jump when $(L,R)= (u,v)$) and similarly, that $\zeta_{[u,v]}$ describes the jump intensity measure for $M_t^-$. 

Here comes a simple observation: Suppose that the pair $(u,v)$ is defined under the infinite measure 
$$ \pi (\mathrm{d}u\,\mathrm{d}v)  = \frac {\mathrm{d}u\,\mathrm{d}v}{(v-u)^2}  1_{\{ u < 0 < v \}}$$
on the set $\R_- \times \R_+$ (it is important for what follows that we restrict ourselves to this set!). 
We can then define $w$ under the measure $\zeta_{[v,u]}$ so that the triple $(u,v,w)$ is defined on the set $\{u<0<v< w \}\cup \{w < u < 0 < v \}$ by the measure with density 
$$\nu (\mathrm{d}u\,\mathrm{d}v\,\mathrm{d}w) = \frac { \mathrm{d}u\, \mathrm{d}v\, \mathrm{d}w } {(w-v) (v-u) (w- u)}. $$ 
We recognize here (a multiple of) the Haar measure on unmarked hyperbolic triangles in $\H$ (see Section \ref{preliminaires}), restricted to those triangles that separate $0$ from $\infty$. Note also that when one sees this triangle without knowing which point is $u$ or which point is $v$, one can recover it immediately: These triangles always have at least one point on $\R_+$ and one point on $\R_-$. If there are two points in $\R_+$ then $u< 0 < v < w$ and if there are two points in $\R_-$, then 
$w<u < 0 < v$. 
 
We can use the same procedure in the other direction. Let us first define $(\alpha, \beta)$ under the same measure $\pi$ and then $\gamma$ under $\zeta_{[\alpha, \beta]}$ (mind that this time, $\gamma \in [\alpha, \beta]$) so that one obtains the triple $(\alpha, \beta , \gamma)$ defined on the set
$\{ \alpha <  \gamma < 0 < \beta \} \cup \{ \alpha < 0 < \gamma < \beta \}$ under the measure with intensity 
$\mathrm{d} \alpha \mathrm{d} \beta \mathrm{d}\gamma /(( \beta-\alpha)(\gamma-\alpha) (\beta -\gamma))$. In this way, we get exactly the same measure as before, and we can also recover from the unmarked triangle which apexes are $\alpha,\beta$ and $\gamma$. 

We have therefore just proved by looking at the properties of the jump measures of the two processes $M$ and $M^-$ that (see for instance \cite {CW} for background on duality for Markov processes):
\begin {proposition}
\label {invmeasure}
The measure $\pi$ is an invariant measure for the Markov transition kernel of $M$ and the dual kernel is that of $M^-$. 
\end {proposition}
In plain words. If for a given positive $t$:
\begin {itemize}
\item
We define $(u,v)$ according to $\pi$, then sample a Markov process $M$ starting from $M_0 = (u,v)$. This defines an infinite measure on quadrilaterals $(u,v,L_t,R_t)$.    
\item 
We define $(u', v')$ according to $\pi$, then sample a Markov process $M^-$ starting from $M_0^- = (u', v')$. This defines an infinite measure on quadrilaterals $(L_t^-, R_t^-, u', v')$. 
\end {itemize}
Then these two measures on quadrilaterals are the same.

\medbreak
Suppose now that $s=\i y$ is some fixed point on the vertical line with $y > 1$.
Let us first sample a triangle $T_\i$ according to $P_\i^\H$ (this is the probability measure obtained by restricting the Haar measure on triangles in $\H$ to those that contain $\i$).
Define its three apexes by $a, l, r$ in such a way that $ l < a < r$. The arc of $T_\i$ that separates $\i$ from infinity is therefore $(l,r)$. Note that one way to sample $T_\i$ is to define $l,r$ under some universal constant $c_{0}$ times  $\pi$, and then $a$ under $\zeta_{[l,r]}$ and to finally restrict the obtained measure on $(l,a,r)$ to those triangles that contain $\i$. 

We now define a triangle $T(s)$ that contains $s$. If $s \in T_\i$, we take $T(s)=T_\i$. 
Otherwise, we sample an accordion $ \mathcal{A}_{(l,r) \to \infty} $. Then, almost surely, $s$ is in one of the triangles of the accordion, that we call $T(s)$. 
Our goal in this paragraph is to show the following reversibility property of the joint law of $(T_\i , T(s))$:

\begin {lemma}
\label {exchange}
If $\varphi : z \mapsto -y  / z$ denotes the M\"obius transformation in $\H$ that interchanges $\i$ and $s$, the law of $(T_\i, T(s))$ is equal to that of $(\varphi ( T(s)), \varphi (T_\i))$. 
\end {lemma}
An equivalent possibly clearer way to phrase this reversibility goes as follows: Define a random triangle $T_s'$ that contains $s$ according to $P_s^\H$, and then define the triangle $T' (\i)$ that contains $\i$, obtained when letting an accordion grow from $T_s'$ towards $0$. Since $( T' ( \i), T_s')$ is distributed as  $(\varphi ( T(s)), \varphi (T_\i))$, the lemma says that 
$( T_\i, T(s))  $ and $(T'(\i), T_s')$ are identically distributed.

\proof
Let us first notice that when restricted to the event that the two triangles are equal (i.e., $s \in T_\i$ for the first pair, and $\i \in T_s'$ for the second one), the laws of 
$( T_\i, T(s))  $ and $(T'(\i), T_s')$ are equal (they are both described via the Haar measure on triangles restricted to the triangles that contain both $\i$ and $s$). We can therefore focus on the event where the two triangles are different.

Suppose now that $T_\i$ (and therefore $(l,r)$) has been sampled and does not contain $s$. When one grows the accordion from the arc $(l,r)$ towards $\infty$, we use the Poisson point process described in previous subsections. Note that we are interested in the law of the first ``time'' (in the Poisson point process parameterization) at which one discovers a triangle (i.e., a jump in our pure jump process) that swallows the point $s$. Classical theory for Poisson point processes (see for instance the ``master formula'' in \cite {RY99}) shows that it is possible to decompose the law of $T(s)$ according to the time at which the jump over $s$ occurs. More precisely, let us grow the accordion from $M_0=(L_0, R_0):= (l,r)$ towards infinity, and let $M_t=(L_t, R_t)$ denote the top boundary arc of the accordion at time $t$.

For each given $t$, we can sample $M_t$ and then define a point $w$ on $\R \setminus [L_t, R_t]$ according to the measure $\zeta_{[R_t, L_t]}$. This is an infinite measure, but the mass $m_t$ of the event that $s$ is in the triangle $(L_t, R_t, w)$ is finite, and furthermore, the definition of the Poisson point process, together with the fact that the triangulation is complete ensures that $E ( \int_0^\infty m_t  \mathrm{d}t ) = 1$ (just because there is a.s. exactly one triangle in the accordion that contains $s$), given the fact that $s \notin T_\i$.
 
Then (conditionally on $(l,r)$ and on the fact that $s$ lies above this arc), the distribution of $T(s)$ is described by 
 \begin{eqnarray*} E \Big( g( T(s)) \,\left| (l,r)\Big) \right. &=&  \int_0^t \mathrm{d}t E_{(L_0,R_0)= (l,r)}  \left(  \int  \zeta_{[R_t, L_t]}(\mathrm{d}w) g(\{L_t, R_t , w\}) 1_{\{ s \in (L_t, R_t, w) \}} \right),  \end{eqnarray*}
for any measurable bounded function $g$. For convenience, we are now going to assume that $g(T) = 0$ as soon as $s \notin T$, as this will enable us to incorporate the indicator function in $g$. As we anyway restrict ourselves to the case where $\i \notin T(s)$, we assume as well that $g(T) = 0$ as soon as $\i \in T$. Similarly, we will consider a measurable function $f$ on the set of unmarked triangles, such that $f(T) = 0$ as soon as $\i \notin T$ or  $s \in T$.  

If we now combine this with the description of the law of $T_\i$, we get that
\begin {eqnarray*}
\lefteqn {E( g(T(s)) f (T_\i )) }
\\
&=& 
\int_0^\infty \mathrm{d}t  c_0 \pi ( \mathrm{d}\alpha  \mathrm{d}\beta)  E_{(L_0,R_0)= (\alpha,\beta)} \left( \int \zeta_{[\alpha,\beta]} ( \mathrm{d}\gamma) \zeta_{[R_t, L_t]}(\mathrm{d}w) f (\{\alpha,\beta,\gamma\}) g (\{L_t, R_t, w\}) \right) .
\end {eqnarray*}
It now suffices to apply the reversibility of our Markov process i.e. Proposition \ref {invmeasure} which implies that this quantity is equal to 
$$
\int_0^\infty \mathrm{d}t  c_0 \pi (\mathrm{d}u \mathrm{d}v)  E_{(L_0^-, R_0^-) = (u, v)}  \left( \int \zeta_{[L_t^-,R_t^-]} (\mathrm{d}\gamma) \zeta_{[v, u]}(\mathrm{d}w) f (\{\gamma, L_t^-,R_t^-\}) g (\{u,v,w\}) \right)
$$ 
and to note that (just in the same way as before), this is precisely equal to $ E ( g ( T_s') f ( T' (\i)))$ which completes the proof.
\endproof

\subsection {End of the proof of Theorem \ref {main}.} \label{section:end} We are now finally ready to conclude the proof of our main result:

\proof 
Let us construct a triangulation of $\TT$ in the unit disk as follows. First sample $T(0)$ according to $P_0$. In each of the three remaining domains $O_1$, $O_2$ and $O_3$ (see Fig.\,\ref{notations}) sample independently an accordion tree started from $(u_2 u_3)$, $(u_3u_1)$ and $(u_1u_2)$ respectively. This clearly defines a random complete triangulation $\TT$. 
Furthermore, this construction ensures directly that $\TT$ is Markovian. What remains to be checked is its M\"obius invariance. 

As the definition also yields invariance under rotations around the origin, it suffices to check that for any given $z_0$ in $\D \setminus \{0\}$, the law of $\TT$ is invariant under the hyperbolic isometry $\varphi_{z_0}$ in $\D$ that interchanges $0$ and $z_0$. 

In other words, define another random triangulation as follows. First sample $T'(z_0)$ according to $P_{z_0}$. Then, in each of the three remaining domains, sample independently an accordion 
tree. This defines a random triangulation $\TT'$ (and by construction, its law is that of $\varphi_{z_0} ( \TT)$ because the image measure of $P_0$ under $\varphi_{z_0}$ is $P_{z_0}$ and by invariance of the accordion tree under M\"obius transformations). Our goal is to prove that $\TT$ and $\TT'$ are identically distributed. 

{
As in the proof of Lemma \ref{twodimmar}, in order to prove that $\TT$ and $\TT'$ are identically distributed, it  suffices to see that for all $Z= \{ 0, z_0, z_1, \ldots, z_n \}$ and for any tree structure $\Gamma$ on $Z$, one has identity in law between $T(Z) 1_{T (Z) \in \mathcal {A} ( Z, \Gamma)}$ and $T'(Z) 1_{T '(Z) \in \mathcal {A} ( Z, \Gamma)}$.

But Lemma \ref {exchange} shows readily that this is certainly the case as soon as $0$ and $z_0$ are neighbors in $\Gamma$. Indeed, we can note that for both triangulations the marginal law of  $(T(0), T(z_0))$ is the same (this is just the image in $ \mathbb{D}$ of Lemma \ref {exchange}), and that (on the event that the remaining points lie in the four outer connected components of the complement of these two triangles) the conditional distribution of $T(Z)$ given these two triangles is also the same by construction.

Suppose now that in $\Gamma$, $0$ is at distance $2$ of $z_0$, i.e. that both $0$ and $z_0$ are neighbors of some $z_l$. Then we can use another random triangulation $\TT''$ that is constructed just as $\TT$ and $\TT'$, except that it is based at $z_l$, i.e., one starts defining $T''(z_l)$ using $P_{z_l}$ etc. 
We have just showed that
\begin {itemize}
\item  $T(Z) 1_{T (Z) \in \mathcal {A} ( Z, \Gamma)}$ and  $T''(Z) 1_{T (Z) \in \mathcal {A} ( Z, \Gamma)}$ are identically distributed.
\item  $T''(Z) 1_{T'' (Z) \in \mathcal {A} ( Z, \Gamma)}$ and  $T'(Z) 1_{T '(Z) \in \mathcal {A} ( Z, \Gamma)}$ are identically distributed.
\end {itemize}
which clearly settles the case when in the case where the graph-distance between $0$ and $z_0$ is $2$.
The very same argument works without any problem to deal with the case where the distance in $\Gamma$ between $0$ and $z_0$ is larger. }
\endproof

\section {The Markovian hyperbolic quadrangulations, pentangulations etc.}
\label {S5}

It is of course natural to wonder if our results are specific to hyperbolic triangulations, or if they do have extensions to random partitions into other hyperbolic polygons. The answer to this question is that, while triangulations are of course in some way special,
 there exist analogs to our Markovian triangulation when one replaces triangles by other polygons. In order to avoid any notational mess, the discussion in the present section will deliberately remain on a rather descriptive with a maybe wordier style, and we will leave mathematical details to the interested reader. 

Let us first focus on the case of \emph{regular polygons}: We say that the sequence of distinct points $x_1, \ldots, x_n$ on the unit circle that is ordered anti-clockwise is a hyperbolic
$n$-gon (and the $n$ hyperbolic lines $(x_1x_2)$ etc. are its boundary edges). We say that it is a \emph{regular} $n$-gon if there exists a M\"obius transformation $\phi$ such that for all $j$ in $\{1, \ldots, n \}$, $\phi ( x_j) = e^{2 \mathrm{i}\pi j/n }$.
 For instance, any (anti-clockwise ordered) triple $x_1, x_2, x_3$ is a regular $3$-gon, but then, only one possible $x_4$  turns $x_1, x_2,x_3, x_4$ into a regular  $4$-gon (we will from now on call regular $4$-gons \emph{hyperbolic squares}). Note that any $n$-gon is obtained by glueing together $(n-2)$ adjacent triangles (and that triangles are all equivalent up to hyperbolic isometry). The $\mu$-hyperbolic area of any $n$-gon is therefore always equal to $n$ (with our normalization of the hyperbolic measure $\mu$).
A hyperbolic square is just the glueing of two ``conformally symmetric'' adjacent triangles and it has $\mu$-area equal to $2$. Note also that an (unmarked) hyperbolic $n$-gon corresponds to $n$ different possible marked $n$-gons (one has to choose which one of the corners is $x_1$). 

It is trivial to extend the definition of the Markovian property to random tilings of $\D$ into $n$-gons. 
The first extension of Theorem \ref {main} goes as follows:
 
\begin {theorem}
For any $n \ge 3$, there exists exactly one (law of a) Markovian M\"obius-invariant complete partition of $\D$ into regular $n$-gons.
\end {theorem}
 
The case $n=3$ is exactly Theorem \ref {main}, whereas when $n=4$, the statement is that there exists a unique Markovian M\"obius-invariant partition of $\D$ into hyperbolic squares. 

\proof[Sketch of the proof.]
The proof of this theorem goes along similar lines as Theorem \ref {main}. Let us focus here on the case $n=4$ (the proof of the other cases is quasi-identical and involves essentially no other idea) and just highlight the main differences with the proof of the case $n=3$. \medskip

\textsc{The origin square.} First of all, note that there is a natural measure $\nu_\circ^{(4)}$  on marked hyperbolic  squares: It is the measure that is obtained from the measure $\nu_\circ$ on marked triangles $(a,b,c)$ by looking at $(a,b,c,d)$ where $d$ is the symmetric image of $b$ with respect to $(ac)$ i.e., the only point $d$ such that $(a,b,c,d)$ is a hyperbolic square. Clearly, this measure is invariant under circular relabeling (i.e., under $(a,b,c,d) \mapsto (b,c,d,a)$), and the marginal measure of any of its four halves (i.e.\,the triangles obtained by dropping one of the four points) is $\nu_\circ$. It follows immediately that if a random tiling $ \mathbf{S}$ of the disk into hyperbolic squares is M\"obius-invariant, then the measure on the marked square $S(0)$ containing the origin (if one marks it by choosing one of its corners uniformly at random among the four) is a multiple of $\nu_\circ^{(4)}$ restricted to those squares that contain the origin (the multiplicative constant being chosen in such a way that the $\nu_\circ^{(4)}$-mass of 
the set of marked squares containing the origin is equal to $1$).  \medskip

\textsc{Uniqueness.} Let us now consider the accordion of squares in the quadrangulation $ \mathbf{S}$ between $0$ and $1$ in $\D$. There is still no problem to define this object; it corresponds to the set of squares of $ \mathbf{S}$ that intersect the segment $[0,1]$. 

\begin{figure}[!ht]
\begin{center}
\includegraphics[width=13.4cm]{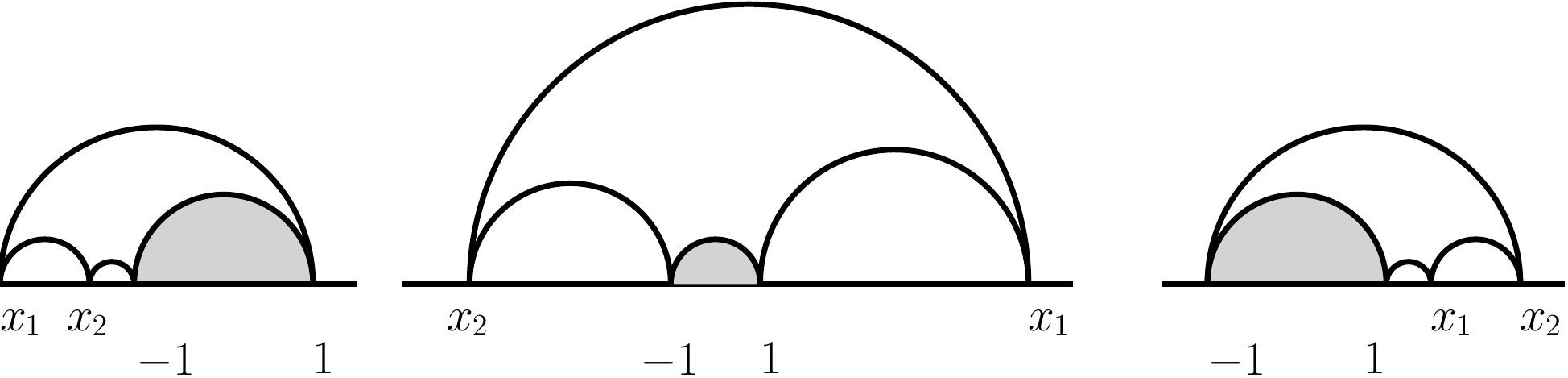}\\
\caption{Sketch of squares of type I$_1$, II and I$_2$. \label{square}}
\end{center}
 \end{figure}

In the uniqueness part of the proof, we assume that $ \mathbf{S}$ is a M\"obius-invariant Markovian decomposition into hyperbolic squares. One then first checks using the Markovian property and M\"obius invariance, that the corresponding half-plane accordion $\mathrm{Acc}_{ \mathbf{S}} ( \mathrm{i}, \infty)$ can be described via a Poisson point process of intensity given by some measure $\rho^{(4)}$ on hyperbolic squares in $\HH$ that have two adjacent corners at exactly $-1$ and $1$.
The measure $\rho^{(4)}$ plays the same role as $\rho$ in Section \ref {S3}; it can be viewed as a measure on the set of pairs of points $(x_1, x_2)$ in $\R \setminus [-1, 1]$ such that $-1, 1, x_1, x_2$ are ordered anti-clockwise on $\partial \HH$, and such that $-1, 1, x_1, x_2$ is always a hyperbolic square, see Fig.\,\ref{square}. The fact that $(-1, 1, x_1, x_2)$ is a hyperbolic square means that
 \begin{eqnarray} \label{*}  \frac {x_2 -1}{x_2+1} &=& 2 \times \frac {x_1 - 1}{x_1+1}.  \end{eqnarray}

Just as for $\rho$ (using M\"obius invariance), one checks that the image of this measure under any M\"obius transformation that preserves $-1$ and $1$ is a multiple of $\rho^{(4)}$. If one now defines $\rho_1^{(4)}$ to be the image measure of $\rho^{(4)}$ under the mapping $(x_1, x_2) \mapsto x_1$, it follows readily (using analogous arguments as in Section \ref {S3} in order to show that $\rho$ is in fact invariant under those transformations) that $\rho^{(4)}_1$ is equal to $\rho$ itself (or to a constant multiple of $\rho$) and that it is equal also to the image of $\rho^{(4)}$ under $(x_1, x_2) \mapsto x_2$. Hence the measure $\rho^{(4)}$ is totally described by \eqref{*} and the fact that $\rho^{(4)}_{1}= \rho$ (up to multiplicative constant).

The rest of the uniqueness part of the proof (the fact that this $\rho^{(4)}$ characterizes the law of the accordion, and that the law of the accordion characterizes the law of $ \mathbf{S}$) follows exactly the same arguments as in the case $n=3$ in Section \ref {S3}.  \medskip

\textsc{Existence.} 
We can decompose the set of squares $(-1, 1, x_1, x_2)$ into three parts: The part I$_2$ where $1< x_1 < x_2$ (and on this part, the only relevant information in order to construct the ``future'' of the accordion is the value of $x_2$, which is distributed just as $\rho$ on $[1, \infty)$), the part I$_1$ where $x_1< x_2< -1$ (and here, the only relevant information for the future of the accordion in $x_1$, which is distributed just as $\rho$ on $(-\infty, -1]$), and the third part that we will call II with $x_2 < -1 < 1 < x_1$ (and here, one needs to know both $x_1$ and $x_2$ to construct the future of the accordion).

{
We then consider a Poisson point process of squares with intensity $\rho^{(4)}$. For a square $S$ in the upper half-plane (with four apexes on the real line), denote its left-most apex by $a^{l} (S)$ and its right-most one $a^r (S)$. This therefore defines a Poisson point process $\{(a^l_{t_i}, a^r_{t_i})\}$ (the square corresponding to $(a^l_t, a^r_t)$ is of type I or II depending on whether one of the values 
$-a^l_t$ and $a^r_t$ is equal to $1$ or not). Then, as in the case of the triangulations, one defines two pure jump processes $(L_{t},R_{t})$ by putting 
 \begin{eqnarray*} (L_t, R_t) &=& \varphi_{t-}^{-1} ((a^l_t, a^r_t)),  \end{eqnarray*}
at all jump-times (where $\varphi_{t-}^{-1}$ is the affine map that maps $(-1, 1)$ onto $(L_{t-}, R_{t-})$). Contrary to the triangulation case, $L$ and $R$ can jump simultaneously (if the corresponding square is of type II).
Note that I${}_1$ and I${}_2$ have infinite $\rho^{(4)}$-mass. The key-observation is that the $\rho^{(4)}$-mass of type-II squares is finite (this can be for instance seen from the fact that this mass is equal to $\rho(\{4 \leq x_{1} \})$).}
Hence, when one constructs the half-plane accordion, the ``times'' (in the Poisson point process) at which one discovers a square where the two sides separating the square from $0$ and from $\infty$ are not adjacent, form a discrete locally finite set. In other words, the accordion will be quite similar to the accordion with triangles except that:
\begin {itemize}
\item For each triangle $(abc)$ in the accordion, one adds the fourth point $d$ in order to turn it into a square, in such a way that $d$ does lie on the same side of the triangle $(abc)$ as neither $0$ nor infinity. 
\item One has squeezed in, in a Poissonian way, a discrete locally finite family of hyperbolic squares where the sides that separate $0$ from infinity are not adjacent.   
\end {itemize}
It then suffices to use the fact that the ``times'' at which one adds those type-II squares is in fact exactly the same (modulo time-reversal) when one looks at the accordion in $\D$ from $-1$ to $1$ or from $1$ to $-1$ (here, one uses the fact that the ``time''-parametrizations of the forward and backward accordions are the same, except for the time-reversal). This will then indeed allow to prove the M\"obius invariance of the random decomposition into squares that is defined in this way.
\endproof

Let us now list further possible extensions:
\begin {itemize}
\item 
The previous proof shows that it is also very easy to construct a tiling of $\D$ into a mixture of triangles and squares. The M\"obius-invariant measure would then be described by a parameter $p \in [0,1]$ (each value of $p$ would correspond to a different distribution on tilings) that is equal to the probability that (in the tiling) the origin is in a triangle (and not in a square). More generally, for any distribution $P$ on $\{3, 4, 5, \cdots, n_0 \}$, one can define a random tiling of $\D$ into regular $n$-gons (with varying $n$'s) in such a way that the probability that the origin lies in a regular $n$-gon is equal to $P(n)$. Conversely, these tilings will be the only complete Markovian tilings of $\D$ into regular $n$-gons for $n \le n_0$. 
\item    
So far, we have been dealing with tilings by regular $n$-gons. When $n=3$, any triangle is a regular $3$-gon, so that this was not a restriction, but for $n \ge 4$, it is. 
It raises the additional question of the existence and characterization of Markovian complete tilings of $\D$ into general $n$-gons. It turns out to be very easy to construct such tilings by non-symmetric tiles. For instance, we could be looking for tiles that are $4$-gons with a prescribed conformal structure i.e. such that the four ordered boundary points $(x_1, x_2, x_3, x_4)$ can be mapped by some M\"obius transformation onto one of the two $4$-gons $(i, -1 , -i, e^{ \mathrm{i} \theta})$ or  $(i, -1 , -i, e^{-\mathrm{i}\theta})$ for some given $\theta \in (0, \pi/2)$ (note that the condition has also to be satisfied also by $(x_2, x_3, x_4, x_1)$). This is then the unique Markovian tiling of the disk into hyperbolic rectangles of prescribed aspect ratio. 
 Loosely speaking, the only main difference with the previously described case of tiling by squares is now that in the accordion, one tosses a fair coin for each $4$-gon
 in order to 
choose between $\theta$ or $-\theta$ (i.e. if one discovers one of its ``long'' sides or one of its ``short'' sides). 

The general statement about Markovian tilings by $n$-gons (of non-necessarily prescribed hyperbolic structure) would then go along the following lines:
Suppose that $J$ is a M\"obius-invariant measure on the set of ordered polygons $P_n$ where 
$$P_n  = \{ (x_1, x_2, \ldots ,x_n ) \in \partial \D, \ \hbox { ordered clockwise}\}$$
 such that 
\begin {enumerate}
\item The image measure of $J$ via the projection $(x_1, \ldots, x_n) \mapsto (x_1, x_2, x_3)$ of $P_n$ onto $P_3$ is a multiple of $\nu_\circ$.
\item The $J$-mass of the set of $n$-gons that contain the origin is equal to one (this is just a matter of normalization). 
\item $J$ is invariant under $(x_1, x_2, \ldots, x_n) \mapsto (x_2, x_3, \ldots, x_n, x_1)$.  
\end {enumerate}
Then, there exists a complete M\"obius-invariant Markovian tiling of $\D$ such that the $n$-gon containing $0$ (if one marks it by choosing $x_1$ uniformly
 at random among the corners of the $n$-gon) is distributed according to the restriction of $J$ to those $n$-gons that contain the origin. Conversely, this construction basically describes all possible complete Markovian M\"obius-invariant tilings by $n$-gons. 

Furthermore, all such measures $J$ can be obtained via a product measure $\nu_\circ \otimes P$, where $\nu_\circ$ chooses the first three points $x_1, x_2, x_3$ and $P$ the hyperbolic position of the other $n-3$ points with respect to $(x_1, x_2, x_3)$. 

We leave out the details of the proofs, as well as the generalizations to tilings into mixtures of $n$-gons for varying $n$'s to the interested reader. 
 \end {itemize}

\section {Concluding remarks}
\label{Scomments}

We conclude the paper with some remarks and open questions.

\paragraph{Hausdorff dimension.}
One property of our Markovian triangulation that we have collected on the way (we safely leave the details to the reader, recall Lemma \ref{box} and Corollary \ref{Kinfini}) is that:
\begin {proposition}
The Hausdorff dimension of the closure of the union of all triangle boundaries in our Markovian hyperbolic triangulation is almost surely equal to $1$.
\end {proposition}
In other words, the ``dimension'' of the Markovian triangulation is not larger than $1$, as opposed to other natural random triangulations that are ``fatter'' (see for instance \cite {Ald94b,CLG10}). 

\medbreak

\paragraph{On completeness.}
For non-complete M\"obius-invariant triangulations, one can still make sense of the definition of the Markovian property in the following way: First note that the probability that $z \in \mathbb{D}$ is in some triangle of the triangulation is equal to some constant $p_{0}$ that does not depend on $z \in \mathbb{D}$ because of M\"obius invariance, and if we assume that $\TT$ is not almost surely empty, $p_{0}$ is strictly positive. Then, we can condition on the event that $T(0)$ is not empty, and define the Markovian property as in Section \ref{markov}.

It is easy to define non-trivial M\"obius-invariant Markovian triangulations that are not complete. Here is an example: Pick $p \in [0,1)$ and consider our random complete triangulation $\TT$ as defined in Theorem \ref{main}. Conditionally on $\TT$, let $( \mathsf{d}_{T})_{T \in \TT}$ be independent Bernoulli variables of parameter $p$ indexed by the triangles of $\TT$, and define
 \begin{eqnarray*} \TT^{(p)} = \left\{ T \in \TT : \mathsf{d}_{T}=1 \right\}. \end{eqnarray*}
In other words, we keep each triangle of $\TT$ with probability $p$. Clearly, $\TT^{(p)}$ is a non-complete M\"obius-invariant Markovian triangulation.
 
It is nevertheless possible to strengthen Theorem \ref {main} replacing the completeness assumption by a density assumption. We recall that a triangulation $ \mathrm{T}$ is dense if the union of the triangles of $ \mathrm{T}$ is dense in $\mathbb{D}$. The statement then becomes:
{\sl There is a unique (law of a) dense M\"obius-invariant triangulation of $\mathbb{D}$ that fulfills the spatial Markov property.}

In other words, any dense M\"obius-invariant Markovian triangulation is in fact complete, see Remark  \ref{dense-complete}. We chose for expository reasons to focus on Theorem \ref {main} and its proof in the present paper, and we therefore do not include the proof of this last statement here, but let us nevertheless give some brief hints to the interested reader:  Most of the analysis goes along similar lines as that of Section \ref {S3}. Density makes it possible to still define properly the processes $(\ell)$,$(r)$ and the jumps $X$. The fact that the triangulation is not complete however a priori allows the possibility to include a drift part in the subordinator-like evolution of $r$ and $\ell$.  But it turns out that there is no M\"obius-invariant way to define a non-zero drift term, as applying a hyperbolic isometry would change the relative speed of the drift near to $r$ and to $\ell$ (which is the usual problem when one tries to define in a M\"obius-invariant way to define processes growing simultaneously at two different points).

\medbreak

\paragraph {Open questions.}
 We conclude with three of the natural open questions: 
\begin {enumerate}
 \item 
It would be nice and enlightening to have an alternative construction of our Markovian M\"obius-invariant triangulation (for instance via an auxiliary Poissonian model, some allocation idea or some statistical physics arguments) that would explain ``directly'' why it exists.
\item
What are the natural discrete models that one could think of, and that would give rise to such Markovian hyperbolic triangulations in the scaling limit? Note that the Markovian property is somehow reminiscent of Omer Angel's exploration of percolation interfaces in random triangulations \cite {Angel} (that however gives rise to a different scaling limit). This could also provide another (heuristic or rigorous) justification to the existence of our triangulation (just as the discrete percolation model ``explains'' the locality and reversibility properties of SLE$_6$ \cite {LSW}). 
\item
What are the natural and M\"obius-invariant dynamics (if they exist) on the set of triangulations that leave our measure invariant? One would for instance like to have a ``continuous'' evolution (in some appropriate topology)  and such that the evolution of triangles that are far away from each other de-correlate fast with this distance. Each ``individual'' triangle could for instance follow some Brownian motion in the three-dimensional Lie group of M\"obius transformations.

\end{enumerate}

\bigbreak

\noindent \textbf{Acknowledgments.} We thank Fr\'ed\'eric Paulin, Pierre Pansu and -- of course -- Itai Benjamini for stimulating discussions. We also thank the referees for their careful reading and useful comments. The authors acknowledge support of the Fondation Cino del Duca and of the Agence Nationale pour la Recherche (projet MAC2).

\null
\bigbreak

\noindent
\begin {tabular}{l p{3mm} l}
D\'epartement de Math\'ematiques et Applications & &Laboratoire de Math\'ematiques 
\\
Ecole Normale Sup\'erieure, 45 rue d'Ulm & &Universit\'e Paris-Sud, B\^at. 425
\\
75230 Paris cedex 05, France & &
91405 Orsay cedex, France
\end {tabular}

\medbreak

\noindent
nicolas.curien@ens.fr

\noindent
wendelin.werner@math.u-psud.fr

 \end{document}